\documentclass[11pt,a4paper,leqno]{article}
\usepackage[T1]{fontenc}
\usepackage[utf8]{inputenc}
\usepackage{lmodern,amsmath,amssymb,amsthm,mathtools}
\usepackage[margin=30mm]{geometry}
\usepackage{microtype}
\usepackage{enumitem,needspace}
\usepackage{titlesec}
\usepackage[colorlinks=true,linkcolor=blue,citecolor=blue,urlcolor=blue]{hyperref}
\hypersetup{pdftitle={Log-Sobolev inequality for the sine-Gordon model in the subcritical regime},
pdfauthor={Alexandru F. Radu}}
\titleformat{\subsection}[runin]{\normalfont\bfseries}{\thesubsection.}{.6em}{}[.]
\titleformat{\subsubsection}[runin]{\normalfont\itshape}{\thesubsubsection.}{.6em}{}[.]
\titlespacing*{\subsection}{0pt}{2.2ex plus 1ex minus .2ex}{.7em}
\titlespacing*{\subsubsection}{0pt}{1.8ex plus .5ex}{.7em}
\numberwithin{equation}{section}
\newtheorem{theorem}{Theorem}[section]
\newtheorem{proposition}[theorem]{Proposition}
\newtheorem{lemma}[theorem]{Lemma}

\theoremstyle{definition}

\theoremstyle{remark}

\DeclareMathOperator{\Hess}{Hess}
\DeclareMathOperator{\Ent}{Ent}
\DeclareMathOperator{\osc}{osc}
\DeclareMathOperator{\Tr}{Tr}
\newcommand{\id}{\mathrm{id}}
\newcommand{\E}{\mathbb E}
\newcommand{\R}{\mathbb R}
\newcommand{\Z}{\mathbb Z}
\newcommand{\N}{\mathbb N}
\newcommand{\eps}{\varepsilon}

\title{Log-Sobolev inequality for the sine-Gordon model\protect\\in the subcritical regime}
\author{Alexandru F. Radu\thanks{Simion Stoilow Institute of Mathematics, Romanian Academy,
Calea Grivitei Street, no.~21, 010702 Bucharest, Romania.
\texttt{alexandru.radu@imar.ro}}}
\date{}
\begin{document}
\maketitle
\begin{abstract}
We prove a logarithmic Sobolev inequality for the massive
two-dimensional sine-Gordon model throughout the subcritical range
$0<\beta<8\pi$. For fixed mass, volume and arbitrary real coupling,
the constant has a positive lower bound independent of the lattice spacing.
When the dimensionless coupling is small, the lower bound is also
uniform in volume and approaches the Gaussian constant as the
coupling tends to zero. This extends the known range $\beta<6\pi$
and gives a rate of entropy decay for the reversible Langevin
dynamics that remains positive as the cutoff is removed.
The proof follows the Polchinski flow to a positive physical scale,
where the Hamiltonian is a bounded perturbation of a uniformly
convex one. For interactions invariant under constant field shifts,
spatial differences bound their curvature by the gradient energy.
For the oscillating part, we compare Gaussian integration with an
average in which the original cosine phase fluctuates independently
of the remaining variables. We estimate the comparison error at an
intermediate scale, then differentiate the remaining
Gaussian density. This preserves both phase damping and spatial
decay without assuming bounds on second derivatives of the incoming
local terms. The resulting strict gain between scales closes the
estimates for every fixed $\beta<8\pi$.
\end{abstract}

\section{Introduction and main result}

We consider the relaxation to equilibrium of the massive sine-Gordon
model in two dimensions as the lattice spacing tends to zero. We prove
a logarithmic Sobolev inequality which remains uniform in this limit
for $0<\beta<8\pi$. Bauerschmidt and Bodineau~\cite[Theorem~1.6]{BB2021}
established the inequality for $\beta<6\pi$ using the Polchinski flow.
To reach the range proposed in~\cite[Section~6.1, p.~262]{BBD2024}, we need
a sharper estimate of the curvature produced by the nonlinear part
of this flow. We average the fluctuations up to a fixed physical
length, compare the remaining Hamiltonian with a convex one, and
transfer its inequality to the original measure.

Gaussian averaging reduces the original cosine but also produces
interactions which depend on field differences. These are unchanged
by a common field shift, so the same decay cannot hold for the whole
potential. We estimate their curvature by the gradient energy and
keep the Gaussian decay for the remaining oscillation. The oscillating
part can generate new neutral terms at each integration, so both
estimates must be preserved.

Let $m,L>0$, $L/\eps\in\N$, and
$\Lambda=(\Z/(L/\eps)\Z)^2$. A site $x$ has physical position $\eps x$.
Write $(f,g)=\sum_xf_xg_x$, let $e_1,e_2$ be the coordinate unit
vectors, and use the nearest-neighbour Laplacian
$\Delta_\Lambda\varphi_x=\sum_{i=1}^2
(\varphi_{x+e_i}+\varphi_{x-e_i}-2\varphi_x)$. For $z\in\R$ and
$0<\beta<8\pi$, put
\begin{equation}\label{eq:unit-lattice}
 A=-\Delta_\Lambda+\eps^2m^2\id,\qquad
 \alpha=2-\frac\beta{4\pi},\qquad
 V_0(\varphi)=2z\eps^\alpha\sum_x\cos(\sqrt\beta\varphi_x).
\end{equation}
The sine-Gordon measure on $\R^\Lambda$ is
\begin{equation}\label{eq:model}
 d\nu_{\eps,L}(\varphi)=Z_{\eps,L}^{-1}
 e^{-\frac12(\varphi,A\varphi)-V_0(\varphi)}\prod_xd\varphi_x.
\end{equation}
In $V_0$, the factor $\eps^\alpha=\eps^2\eps^{-\beta/(4\pi)}$
combines the area element with the compensation for logarithmic
Gaussian variance, as in~\cite[(1.15)--(1.16)]{BB2021}.
Throughout, $\partial_x,\nabla$ and $\Hess$ differentiate field values.
We use the Euclidean norm $\|v\|_2=(\sum_x|v_x|^2)^{1/2}$ for vectors
and its induced matrix norm $\|\cdot\|_{\mathrm{op}}$.
For a field-dependent vector or matrix, $\|\cdot\|_{\infty,2}$ or
$\|\cdot\|_{\infty,\mathrm{op}}$ also takes the supremum over fields.
In particular,
\[
 \|\Hess F\|_{\infty,\mathrm{op}}
 =\sup_{\varphi\in\R^\Lambda}\sup_{\|v\|_2=1}
   \|\Hess F(\varphi)v\|_2,
\]
where the first supremum is over all field configurations and the second
over all Euclidean unit vectors in $\R^\Lambda$.
For scalar functions, $\|f\|_\infty=\sup_\varphi|f(\varphi)|$.
The letters $C,C_\beta,C_M$ denote scalar bounds; a generic Gaussian
covariance will be denoted by $\Sigma$.

For nonnegative $F$, define
$\Ent_\nu F=\E_\nu(F\log F)-\E_\nu F\log\E_\nu F$ and the physically
normalized Dirichlet form
\begin{equation}\label{eq:entropy-energy}
 \mathbf D_{\eps,L}(F)=\eps^{-2}\sum_x\E_{\nu_{\eps,L}}
                         [(\partial_xF)^2].
\end{equation}
Let $\gamma_{\eps,L}$ be the largest constant for which
\begin{equation}\label{eq:lsi}
 \Ent_{\nu_{\eps,L}}F
 \le\frac2{\gamma_{\eps,L}}\mathbf D_{\eps,L}(\sqrt F)
\end{equation}
holds for every smooth positive $F$. The factor $\eps^{-2}$ gives
$\gamma_{\eps,L}=m^2$ in the Gaussian case $z=0$.
By the approximation used in~\cite[Section~2.2, p.~2075]{BB2021}, it suffices
to prove~\eqref{eq:lsi} for smooth $F$ bounded above and away from zero.

\begin{theorem}\label{thm:main}
Let $0<\beta<8\pi$. For every fixed $m,L>0$ and $z\in\R$,
\begin{equation}\label{eq:main-fixed-volume}
 \inf_{\eps>0:\,L/\eps\in\N}\gamma_{\eps,L}>0.
\end{equation}
There are $c_\beta,C_\beta>0$ such that, if
$|z|m^{-\alpha}\le c_\beta$, then
\begin{equation}\label{eq:main-small-coupling}
 \gamma_{\eps,L}\ge m^2e^{-C_\beta|z|m^{-\alpha}}
                  \ge m^2-C_\beta|z|m^{\beta/(4\pi)}.
\end{equation}
The constants in this assertion are independent of $L,m,z,\eps$.
\end{theorem}

The two assertions use different stopping lengths. Since $\alpha>0$,
the dimensionless coupling $|z|r^\alpha$ is small at a sufficiently
short physical length for any fixed $z$, independently of the cutoff. This gives
\eqref{eq:main-fixed-volume}. Under the small-coupling assumption,
we can stop at a fixed fraction of the mass length, or of the side
length in a smaller torus, giving~\eqref{eq:main-small-coupling}.

For the relaxation interpretation, let $f_s$ be the density relative
to $\nu_{\eps,L}$ at physical time $s$ for the reversible Langevin
dynamics with Dirichlet form $\mathbf D_{\eps,L}$.
The entropy dissipation identity~\cite[Proposition~2.2]{BBD2024} gives
\[
 \frac{d}{ds}\Ent_{\nu_{\eps,L}}f_s
 =-4\mathbf D_{\eps,L}(\sqrt{f_s})
 \le-2\gamma_{\eps,L}\Ent_{\nu_{\eps,L}}f_s.
\]
The entropy therefore decays exponentially at a rate independent of
the lattice resolution. The same constant controls fluctuations. By
the spectral-gap consequence~\cite[Proposition~2.5]{BBD2024}, for a
smooth test function $f$ on the physical torus,
\[
 \operatorname{Var}_{\nu_{\eps,L}}
       \left(\eps^2\sum_xf(\eps x)\varphi_x\right)
 \le\gamma_{\eps,L}^{-1}\eps^2\sum_x f(\eps x)^2.
\]
The Riemann sum on the right remains bounded as $\eps\to0$.
The logarithmic Sobolev inequality also gives Gaussian concentration
uniformly in the cutoff~\cite[Proposition~5.4.1 and (5.4.2)]{BGL2014}.

\subsection{The estimates needed along the flow}
In physical units, the microscopic curvature can be as low as
$m^2-2\beta|z|\eps^{-\beta/(4\pi)}$, and
$\sup V_0-\inf V_0\le4|z|L^2\eps^{-\beta/(4\pi)}$.
To obtain cutoff-uniform bounds from the convexity and
bounded-perturbation criteria~\cite{BE1985,HoSt1987}, we apply them
after smoothing. We must therefore control the renormalized potential
at every background field, including fields far from equilibrium.

Already the product of two cosines contains
$\cos(\sqrt\beta(\varphi_x-\varphi_y))$, which is unchanged when
the same constant is added to every field value. We call interactions
with this invariance neutral. Differentiating twice in a field direction
$f$ shows the gain that this invariance can provide,
\begin{equation}\label{eq:neutral-pair}
 \left.\frac{d^2}{ds^2}
 \cos\bigl(\sqrt\beta[\varphi_x+sf_x-\varphi_y-sf_y]\bigr)
 \right|_{s=0}
 =-\beta\cos\bigl(\sqrt\beta(\varphi_x-\varphi_y)\bigr)(f_x-f_y)^2.
\end{equation}
For the local power counting, take $f$ to be a heat-smoothed direction
and measure the separation $r$ in units of the smoothing length.
Each difference then supplies one power of $r$. Together with the
pair coefficient $r^{-\beta/(2\pi)}$ and the radial area element,
this gives the integral
\begin{equation}\label{eq:threshold-integral}
 \int_0^1r^{3-\beta/(2\pi)}\,dr
 =\frac1{4-\beta/(2\pi)}<\infty
 \quad\Longleftrightarrow\quad\beta<8\pi.
\end{equation}
Bauerschmidt and Bodineau already prove this pair bound
in~\cite[Proposition~3.6]{BB2021}. Their higher-order estimate uses
a first spatial moment, which restricts it to $\beta<6\pi$
\cite[Lemmas~3.11 and A.4, Section~3.7]{BB2021}. The pair calculation
therefore identifies the expected range. To obtain a logarithmic
Sobolev inequality in this range, we must control the curvature
of the exact interaction after all orders have been combined.

We keep a spatially weighted bound on the full Hessian. This bound
ensures convexity of the conditional Hamiltonians used in the next
integration and also controls the neutral part. Indeed, the Hessian
of a neutral interaction has zero row and column sums, so its
quadratic form can be written in differences of the test vector, as
in~\eqref{eq:neutral-pair}. Along the flow we apply this identity to
heat kernels. At the stopping scale, spatial translation invariance
and a shortest-path estimate bound it by the gradient energy,
as shown in~\eqref{ap:neutral-form-bound}. This contribution changes
the coefficient of the gradient energy but leaves the Gaussian mass
unchanged.

For the oscillating part, we compare Gaussian averaging with an
average in which the original cosine phase fluctuates independently
of the other variables. The latter has the exact cosine decay.
Their difference involves derivatives with this phase held fixed.
These derivatives vanish for the original cosine, and we control
them by a quadratic bound in the coupling. Their sum over sites is
zero. Subtracting the central heat-kernel value then produces a
second spatial difference, by evenness of the kernel. The resulting
time integral, computed in~\eqref{ap:critical-integral}, is finite
throughout $\beta<8\pi$.

The main difficulty is to preserve this first-derivative bound.
Applying the same comparison to a derivative would require second
derivatives of the local term. We make the comparison at the midpoint,
where the error decays as a power of the scale ratio. Differentiating
the remaining Gaussian density gives a derivative bound with this
decay. The incoming first derivatives give a second bound with
stronger spatial decay. Interpolating the two derivative estimates
retains both forms of decay. When the length is multiplied by $M$,
the coefficient of the incoming
quadratic bound grows at most as $M^{3\alpha/2}$, whereas the allowed
quadratic bound grows as $M^{2\alpha}$. This strict gain, made explicit
in Lemma~\ref{lem:local-periodic-derivatives}, closes the induction
for every $\alpha>0$.

Renormalization-group analyses of the sine-Gordon measure
include~\cite{BGN1982,NRS1986,DH1993,DH2000}. Further constructions
using flows below $6\pi$ appear in~\cite{BH2022,GM2026,Vihko2026}.
Pelai\v c~\cite[Proposition~1.1 and p.~2, footnote~2]{Pelaic2025} constructs the measure throughout the
subcritical range for arbitrary real coupling at finite volume.
Meyer~\cite[Theorems~1.1--1.2 and Proposition~1.3]{Meyer2026} obtains order-by-order approximate potentials
throughout the subcritical range, with polynomial field-gradient
dependence, and a construction by forward--backward stochastic
differential equations below $(6/7)8\pi$. The parameter denoted
$\beta^2$ there is our $\beta$. These field-dependent estimates serve
the construction of the measure. Our curvature criterion requires
a Hessian bound uniform over all background fields. To apply it
through an approximate potential, one would also have to control
the Hessian of the remainder with the same uniformity.

For the dynamics, the cutoff-uniform relaxation estimate complements
the well-posedness theory obtained by stochastic quantization.
Local well-posedness in the full subcritical range is proved
in~\cite[Theorem~1.1]{CHS2018}, and global well-posedness below $6\pi$
in~\cite[Theorem~1.1]{BC2026}. Finite-volume a priori bounds throughout the
subcritical range are obtained in~\cite[Theorem~1.1 and Remark~1.2]{ES2026}.

Section~\ref{sec:induction} deduces the theorem from estimates on the
renormalized potential. Section~\ref{sec:propagation-proofs} establishes
them by iteration.

\section{Reduction to estimates on the renormalized potential}
\label{sec:induction}\label{sec:continuation}

We integrate Gaussian fluctuations in order of their diffusion length,
\begin{equation}\label{eq:heat-flow}
 Q_t=e^{-tA/2},\qquad C_t=\int_0^te^{-sA}\,ds,\qquad
 e^{-V_t(\varphi)}=\mathbf E_{C_t}e^{-V_0(\varphi+\zeta)}.
\end{equation}
Here $\mathbf E_\Sigma$ averages over a centred Gaussian vector $\zeta$
with covariance $\Sigma$, while $\varphi$ is fixed. The heat kernel spreads
over a distance of order $\sqrt{1+t}$ in lattice units, which will be
our length scale. At fixed cutoff the potentials are smooth with bounded
derivatives.
The convolution solves the Polchinski equation~\cite[(18)]{Polchinski1984},
in the normalization of~\cite[Proposition~3.4]{BBD2024},
\begin{equation}\label{eq:polchinski}
 \partial_tV_t=\tfrac12\Tr(\dot C_t\Hess V_t)
                      -\tfrac12(\nabla V_t,\dot C_t\nabla V_t).
\end{equation}
We iterate the Gaussian formula as the finite-time evolution of this
equation.
The Hessian which enters the logarithmic Sobolev criterion is weighted
by $\dot C_t$. Since
$(u,v)_{\dot C_t}=(u,\dot C_tv)=(Q_tu,Q_tv)$, as in~\cite[(1.2), (1.11)]{BB2021},
we can express it in the ordinary Euclidean coordinates as
$\Hess[V_t(Q_t\varphi)]=Q_t\Hess V_t(Q_t\varphi)Q_t$.
The same choice lets us repeat Gaussian integration without changing
its form. Indeed,
\[
 Q_tC_rQ_t=\int_0^r e^{-(t+s)A}\,ds=C_{t+r}-C_t,
 \qquad Q_tQ_r=Q_{t+r}.
\]
Adding these independent Gaussian fluctuations therefore gives
\begin{equation}\label{eq:convolution}
 e^{-V_{t+r}(Q_{t+r}\varphi)}
 =\mathbf E_{C_r}e^{-V_t(Q_t(Q_r\varphi+\zeta))}.
\end{equation}

We will need the following estimates up to a fixed physical length. Write
$(Dw)_{x,k}=w_{x+e_k}-w_x$, so that $D^*D=-\Delta_\Lambda$, and
$\osc f=\sup f-\inf f$.

\begin{proposition}\label{prop:finite-scale}
For each $0<\beta<8\pi$ there are $c_\beta,C_\beta>0$ with the
following property. Let $\ell\ge1$, $t=\ell^2-1$, and suppose
\[
 \eps\ell\le c_\beta\min\{m^{-1},L\},\qquad
 |z|(\eps\ell)^\alpha\le c_\beta.
\]
Then
\begin{equation}\label{eq:finite-full-hessian}
 \|\Hess[V_t(Q_t\cdot)]\|_{\infty,\mathrm{op}}
 \le C_\beta|z|\eps^2(\eps\ell)^{-\beta/(4\pi)}.
\end{equation}
Moreover, $V_t(Q_t\varphi)=N(\varphi)+R(\varphi)$, where $N$ is
invariant under constant shifts and, for all $\varphi,w$,
\begin{align}
 |(w,\Hess N(\varphi)w)|
 &\le C_\beta|z|(\eps\ell)^\alpha\|Dw\|_2^2,
 \label{eq:finite-neutral-form}\\
 \osc R&\le C_\beta|z|L^2(\eps\ell)^{-\beta/(4\pi)}.
 \label{eq:finite-remainder}
\end{align}
The constants are independent of $\ell,\eps,L,m,z$.
\end{proposition}

For small $|z|(\eps\ell)^\alpha$, the bound~\eqref{eq:finite-neutral-form}
lets us absorb the neutral curvature into the gradient energy.
The bounded part is then treated by~\eqref{eq:finite-remainder}.
In a large torus its oscillation grows with the volume.
Under the small-coupling assumption
we use~\eqref{eq:finite-full-hessian} instead, to keep the whole
Hamiltonian convex.

\begin{proof}[Proof of Theorem~\ref{thm:main} assuming Proposition~\ref{prop:finite-scale}]
For $z=0$, the assertion is the Gaussian inequality~\cite{Gross1975}.
For $z\ne0$, we use the finite-time transport criterion of
\cite[Theorem~5.4 and (5.26)]{BBD2024}, following~\cite{Shenfeld2026}.
If $\mu_t$ is the measure with density proportional to
$e^{-\frac12(\varphi,A\varphi)-V_t(Q_t\varphi)}$, and
$\gamma(\mu_t)$ is its Euclidean logarithmic Sobolev constant, the criterion gives
\begin{equation}\label{eq:finite-criterion}
 \gamma_{\eps,L}\ge\eps^{-2}\gamma(\mu_t)
 \exp\left[-\int_0^t\|\Hess[V_s(Q_s\cdot)]\|_{\infty,\mathrm{op}}\,ds\right].
\end{equation}
These are the rescaled measures in~\cite[(5.28)--(5.31)]{BBD2024}.
The factor $\eps^{-2}$ converts the Euclidean inequality to
\eqref{eq:entropy-energy}.

Choose a physical length $r_*>0$ such that
\[
 r_*\le c_\beta\min\{m^{-1},L\},\qquad
 |z|r_*^\alpha\le c_\beta,
\]
reducing $c_\beta$ if necessary. This choice is independent of the
cutoff. Set
\[
 t_*=\max\{(r_*/\eps)^2-1,0\}.
\]
For $\eps\le r_*$, the assumptions of
Proposition~\ref{prop:finite-scale} hold up to $t_*$, since both the
physical length and its coupling increase with time. Use its
decomposition $V_{t_*}(Q_{t_*}\cdot)=N+R$. The loss in the transfer is
\begin{equation}\label{eq:stopped-curvature}
 \begin{aligned}
 \int_0^{t_*}\|\Hess[V_s(Q_s\cdot)]\|_{\infty,\mathrm{op}}\,ds
 &\le C_\beta|z|\eps^\alpha
       \int_0^{t_*}(1+s)^{-\beta/(8\pi)}\,ds\\
 &=\frac{2C_\beta}{\alpha}|z|\eps^\alpha
       \bigl((1+t_*)^{\alpha/2}-1\bigr)
 \le C_\beta|z|r_*^\alpha.
 \end{aligned}
\end{equation}
If $\eps>r_*$, the integral is zero. In this case put $N=0$ and
$R=V_0$. The bare estimates give
\[
 \begin{aligned}
 \osc V_0&\le4|z|L^2\eps^{-\beta/(4\pi)}
       \le4|z|L^2r_*^{-\beta/(4\pi)},\\
 \|\Hess V_0\|_{\infty,\mathrm{op}}
       &\le2\beta|z|\eps^2r_*^{-\beta/(4\pi)}.
 \end{aligned}
\]
Thus, at every cutoff, $V_{t_*}(Q_{t_*}\cdot)=N+R$ with the same
stopping estimates.
For the part containing $N$, the neutral form bound gives
\[
 A+\Hess N
 \ge\bigl(1-C_\beta|z|r_*^\alpha\bigr)D^*D
          +\eps^2m^2\id
 \ge\eps^2m^2\id,
\]
and the oscillation and full Hessian bounds give
\begin{equation}\label{eq:relative-hessian}
 \begin{aligned}
 \osc R&\le C_\beta|z|L^2r_*^{-\beta/(4\pi)},\\
 \|\Hess[V_{t_*}(Q_{t_*}\cdot)]\|_{\infty,\mathrm{op}}
 &\le C_\beta|z|\eps^2r_*^{-\beta/(4\pi)}.
 \end{aligned}
\end{equation}
Bakry--\'Emery~\cite[Theorem~2.8]{BBD2024} and the Holley--Stroock
perturbation bound~\cite[Proposition~5.1.6]{BGL2014} consequently give
$\gamma(\mu_{t_*})\ge\eps^2m^2e^{-\osc R}$.
Together with~\eqref{eq:finite-criterion}, this yields
\[
 \gamma_{\eps,L}\ge m^2\exp\left[-C_\beta|z|
       \bigl(r_*^\alpha+L^2r_*^{-\beta/(4\pi)}\bigr)\right]>0.
\]

For the volume-uniform assertion, first suppose $Lm\ge1$ and take
$r_*=c_1/m$, with $c_1>0$ fixed and small. The assumption
$|z|m^{-\alpha}\le c_\beta$ makes this choice admissible.
The full Hessian bound in~\eqref{eq:relative-hessian} now gives
\[
 A+\Hess[V_{t_*}(Q_{t_*}\cdot)]
 \ge\eps^2m^2\bigl(1-C_\beta|z|m^{-\alpha}\bigr)\id.
\]
Bakry--\'Emery applies to the whole measure $\mu_{t_*}$, and
\eqref{eq:finite-criterion} gives
\[
 \gamma_{\eps,L}
 \ge m^2\bigl(1-C_\beta|z|m^{-\alpha}\bigr)
           e^{-C_\beta|z|m^{-\alpha}}
 \ge m^2e^{-C_\beta|z|m^{-\alpha}},
\]
after changing $c_\beta,C_\beta$.
If $Lm<1$, take $r_*=c_1L$. The fixed-volume exponent now satisfies
\[
 r_*^\alpha+L^2r_*^{-\beta/(4\pi)}
 =(c_1^\alpha+c_1^{-\beta/(4\pi)})L^\alpha
 \le C_\beta m^{-\alpha},
\]
which gives the same bound.

Finally $e^{-u}\ge1-u$ yields
\eqref{eq:main-small-coupling}.
\end{proof}

\section{Estimates for the renormalized potential}
\label{sec:propagation-proofs}

We now prove Proposition~\ref{prop:finite-scale} by repeating the
Gaussian integration~\eqref{eq:convolution}. The full Hessian bound
will ensure that each conditional Hamiltonian remains uniformly
convex.

\subsection{The decomposition of the potential}
\label{cr:section}\label{ap:section}\label{sec:charged-derivatives}

Fix two length scales $\ell\ge1$ and $\ell'=M\ell$, with $M\ge2$, and put
\begin{equation}\label{ap:coordinates}
 t=\ell^2-1,\qquad T=\ell'^2-\ell^2,\qquad
 Q_tC_TQ_t=C_{t+T}-C_t,\qquad Q_tQ_T=Q_{t+T}.
\end{equation}
Assume $\eps\ell'\le c_0\min\{m^{-1},L\}$, the heat regime of
Lemma~\ref{cr:heat-bounds}, and write $d(i,j)$ for periodic
nearest-neighbour distance. Fix $\mu>0$, depending only on $\beta$,
so that $\max\{32,24/\alpha\}\mu\le\rho_0$.

At length $\ell$, write
\begin{equation}\label{ap:decomposition}
 V_t(Q_t\varphi)=N(\varphi)+R(\varphi),\qquad
 R=\ell^{-2}\sum_xg_x.
\end{equation}
Initially $N=0$ and $g_x=2z\eps^\alpha\cos(\sqrt\beta\varphi_x)$.
The function $g_x$ will change under integration, while the index $x$
continues to refer to this original cosine. Its smoothed phase is
$\sqrt\beta(Q_t\varphi)_x$, which changes under a common field shift
at rate
\[
 c_t=\sqrt\beta\sum_i Q_t(i,x)=\sqrt\beta e^{-\eps^2m^2t/2},\qquad
 \text{period}=2\pi/c_t.
\]
We require $N$ to be neutral and each real $g_x$ to have this period
and zero period average.

The logarithm couples the local terms. To retain the sum over $x$,
we differentiate the logarithm with respect to a parameter multiplying $R$.
For $0\le\lambda\le1$, weight the Gaussian law of $Q_T\varphi'+\zeta$,
with covariance $C_T$, by $e^{-(N+\lambda R)}$ and normalize it.
Denote this expectation by $\mathbf E_\lambda$ and put
\[
 G_x(\varphi')=M^2\int_0^1\mathbf E_\lambda g_x\,d\lambda.
\]
Here $M^2=\ell'^2/\ell^2$ changes the area normalization, and
integration in $\lambda$ gives the exact decomposition,
\begin{equation}\label{fc:first-coupling-identity}
 \begin{aligned}
 &-\log\mathbf E_{C_T}e^{-(N+R)(Q_T\varphi'+\zeta)}
       +\log\mathbf E_{C_T}e^{-N(Q_T\varphi'+\zeta)}\\
 &=\int_0^1\mathbf E_\lambda R\,d\lambda
       =\ell^{-2}\sum_x\int_0^1\mathbf E_\lambda g_x\,d\lambda
       =\ell'^{-2}\sum_xG_x.
 \end{aligned}
\end{equation}
The periods match because
\[
 Q_T(\varphi'+2\pi c_{t+T}^{-1}\mathbf1)
 =Q_T\varphi'+2\pi c_t^{-1}\mathbf1.
\]
Thus $G_x$ has period $2\pi/c_{t+T}$. Its period average,
which is invariant under every constant shift, is
\[
 \Pi f(\varphi')=\frac1{2\pi}\int_0^{2\pi}
        f(\varphi'+\theta c_{t+T}^{-1}\mathbf1)\,d\theta.
\]

We put this average into the neutral interaction and retain the
zero-average part in the local terms,
\begin{equation}\label{eq:next-scale-decomposition}
 \begin{aligned}
 N'(\varphi')&=-\log\mathbf E_{C_T}e^{-N(Q_T\varphi'+\zeta)}
                   +\ell'^{-2}\sum_x\Pi G_x(\varphi'),\\
 g'_x(\varphi')&=G_x(\varphi')-\Pi G_x(\varphi').
 \end{aligned}
\end{equation}
Then~\eqref{eq:convolution} and~\eqref{fc:first-coupling-identity} give
$V_{t+T}(Q_{t+T}\varphi')=N'(\varphi')+\ell'^{-2}\sum_xg'_x(\varphi')$.
The first term in $N'$ is neutral, and $\Pi g'_x=0$. Hence
$N'=\Pi[V_{t+T}(Q_{t+T}\cdot)]$, an identity we shall use to bound
its derivatives by those of the full potential.

The Gaussian average of the original cosine determines the leading
local coefficient. We keep this coefficient exact and use the
dimensionless coupling to bound the nonlinear terms,
\begin{equation}\label{eq:running-coupling}
 z_\ell=z\eps^\alpha e^{-\beta C_t(0,0)/2},\qquad
 \zeta_\ell=|z|\eps^\alpha\ell^\alpha,\qquad
 \ell^2|z_\ell|\le C_\beta\zeta_\ell.
\end{equation}
The inequality follows from~\eqref{cr:diagonal-log}, and the exact
Gaussian average is
\[
 \mathbf E_{C_t}\bigl[2z\eps^\alpha
       \cos(\sqrt\beta(Q_t\varphi+\zeta)_x)\bigr]
 =2z\eps^\alpha e^{-\beta C_t(0,0)/2}
       \cos(\sqrt\beta(Q_t\varphi)_x).
\]
The factor $\ell^{-2}$ in~\eqref{ap:decomposition} makes the leading
local coefficient $2\ell^2z_\ell$, of size $\zeta_\ell$.

To measure the dependence of $g_x$ on variables other than its
original phase, subtract from a coordinate derivative
the common shift producing the same phase change,
\begin{equation}\label{ap:periodic-directions}
 \Theta_t=c_t^{-1}\sum_j\partial_j,\qquad
 \partial_i^{(x)}=\partial_i-\sqrt\beta Q_t(i,x)\Theta_t.
\end{equation}
Indeed, $\Theta_t\cos(\sqrt\beta(Q_t\varphi)_x)
=-\sin(\sqrt\beta(Q_t\varphi)_x)$, whereas
$\partial_i^{(x)}\cos(\sqrt\beta(Q_t\varphi)_x)=0$.
Thus $\Theta_t$ differentiates the angle, whereas $\partial_i^{(x)}$
holds the phase fixed and detects dependence on field differences.

For $g$ with period $2\pi/c_t$ under common shifts,
$\theta\mapsto g(\varphi+\theta c_t^{-1}\mathbf1)$ is $2\pi$-periodic
at each fixed $\varphi$.
Thus $e^{a\Theta_t^2}$ is the heat semigroup on the circle in this
angle, with auxiliary time $a$. On functions of zero shift average,
the supremum-norm estimate is
$\|e^{a\Theta_t^2}g\|_\infty\le Ce^{-a}\|g\|_\infty$.
For the induction we need decay with coefficient one, to retain the
exact Gaussian cosine coefficient while allowing only quadratic errors.
We therefore use the standard equivalent-norm construction for bounded
semigroups~\cite[Lemma~II.3.10]{EngelNagel2000}, applied to
$e^a e^{a\Theta_t^2}$, and set
\begin{equation}\label{ap:circle-norm}
 \begin{gathered}
 \|g\|_* =\sup_{a\ge0}e^a\|e^{a\Theta_t^2}g\|_\infty,\\
 e^{a\Theta_t^2}g(\varphi)
 =\frac1{\sqrt{4\pi a}}\int_\R e^{-\theta^2/(4a)}
       g(\varphi+\theta c_t^{-1}\mathbf1)\,d\theta,\quad a>0.
 \end{gathered}
\end{equation}
For $a\le1$ use supremum-norm contraction, and for $a\ge1$ use
$\sum_{n\ne0}e^{-n^2a}\le Ce^{-a}$. These bounds give norm
equivalence, while the semigroup property gives the required decay,
\[
 \|g\|_\infty\le\|g\|_*\le C\|g\|_\infty,\qquad
 \|e^{b\Theta_t^2}g\|_*
 =e^{-b}\sup_{a\ge b}e^a\|e^{a\Theta_t^2}g\|_\infty
 \le e^{-b}\|g\|_*.
\]
In particular, $\|\cos(\sqrt\beta(Q_t\varphi)_x)\|_*=1$.
Constant-coefficient field derivatives commute with this semigroup, so
the same norm applies to the derivatives of $g_x$. For functions of
$\varphi'$, both derivatives and the norm use $t+T$
in place of $t$. Inside $\mathbf E_\lambda$ they use $t$.

We ask for linear bounds on the full Hessian and the derivative
at $x$. The local norm may exceed the Gaussian cosine coefficient by
a quadratic term, and the phase-fixed derivative is itself quadratic,
\begin{equation}\label{ap:full-hessian-row}
 \sup_i\sum_j e^{\mu d(i,j)/(2\ell)}
       \sup_\varphi|\partial_i\partial_j(N+R)|
 \le H\zeta_\ell\ell^{-2},
\end{equation}
and, at every field,
\begin{equation}\label{ap:norms}
 \begin{aligned}
 \|g_x\|_*&\le2\ell^2|z_\ell|+B\zeta_\ell^2,\\
 \|\partial_i^{(x)}g_x\|_*
 &\le B\zeta_\ell^2\ell^{-2}e^{-6\mu d(i,x)/\ell},\\
 |\partial_xg_x|&\le K\zeta_\ell\ell^{-2}.
 \end{aligned}
\end{equation}
We shall fix $K,H\ge1$ before $B$.
Constants denoted by $C$ depend only on $\beta,\mu$. Constants
$C_M$ may also depend on $M$ and the fixed $K,H$, but not on $B$.
We write $C_{M,K}$ when independence of $H$ is needed. All are
independent of cutoff, volume and field.

\begin{proposition}\label{prop:outgoing-decomposition}
For each $0<\beta<8\pi$, one can choose $M\ge2$, $K,H,B\ge1$ and
$\zeta_*>0$, depending only on $\beta$, so that Gaussian
integration~\eqref{eq:next-scale-decomposition} preserves the full
Hessian bound~\eqref{ap:full-hessian-row} and the local value and
first-derivative bounds~\eqref{ap:norms}, with the same constants
$H,K,B$, from length $\ell$ to $\ell'=M\ell$.
This holds for real $C^2$ functions $N,g_x$, with $N$ invariant under constant shifts
and each $g_x$ periodic under common shifts with period $2\pi/c_t$
and zero period average, provided
\[
 \eps\ell'\le c_0\min\{m^{-1},L\},\qquad
 \zeta_{\ell'}\le\zeta_*.
\]
At the new scale, the derivatives and norm use $t+T=\ell'^2-1$.
\end{proposition}

The proof is completed in Subsection~\ref{sec:completion}, where we
choose the constants. We first note that the derivative at $x$
controls the angular derivative, since $Q_t(x,x)\ge c_\beta\ell^{-2}$.
Indeed, by~\eqref{ap:periodic-directions},
\[
 \Theta_tg_x=\frac{\partial_xg_x-\partial_x^{(x)}g_x}
                       {\sqrt\beta Q_t(x,x)},\qquad
 \|\Theta_tg_x\|_\infty\le C_\beta(K+1)\zeta_\ell
\]
under $B\zeta_\ell\le1$. Since
$\partial_i g_x=\partial_i^{(x)}g_x+\sqrt\beta Q_t(i,x)\Theta_tg_x$,
the pointwise heat-kernel bound gives
\begin{equation}\label{ap:ordinary-gradient}
 \begin{aligned}
 \|\partial_i g_x\|_*
 &\le C\bigl(\|\partial_i^{(x)}g_x\|_\infty
       +\sqrt\beta Q_t(i,x)\|\Theta_tg_x\|_\infty\bigr)\\
 &\le C_\beta(K+1)\zeta_\ell\ell^{-2}e^{-6\mu d(i,x)/\ell}.
 \end{aligned}
\end{equation}
We used equivalence of the circle and supremum norms only on these
derivatives, all of which have zero shift average.

The same Hessian bound applies throughout the interpolation
$N+\lambda R$. To see this, first use the period average to write
\[
 \partial_i\partial_jN(\varphi)
 =\frac1{2\pi}\int_0^{2\pi}
       \partial_i\partial_j(N+R)(\varphi+\theta c_t^{-1}\mathbf1)\,d\theta.
\]
Taking the field supremum inside the integral gives the same weighted
row bound for $N$. The triangle inequality then gives it for
$N+\lambda R=(1-\lambda)N+\lambda(N+R)$, $0\le\lambda\le1$.
In particular, if $\rho/\ell'\le\mu/(8\ell)$,
\begin{equation}\label{ap:source-estimates}
 \sup_i\sum_j e^{\rho d(i,j)/\ell'}
          \sup_\varphi|\partial_i\partial_j(N+\lambda R)|
 \le H\zeta_\ell\ell^{-2}.
\end{equation}
Below we impose $B\zeta_\ell\le1$ and $C_MH\zeta_\ell\le1/8$,
fixing $C_M$ after $M,K,H$ and before $B$ to cover the conditional
estimates. The final choice of $\zeta_*$ will ensure both conditions.

\subsection{Conditional integration}
We shall estimate the covariance terms produced by differentiating a
Gaussian convolution. For a positive definite covariance matrix
$\Sigma$, define
\begin{equation}\label{cr:reference-definitions}
 T_{F,\Sigma}f(\varphi)
 =\frac{\mathbf E_\Sigma[e^{-F(\varphi+\zeta)}f(\varphi+\zeta)]}
        {\mathbf E_\Sigma e^{-F(\varphi+\zeta)}}.
\end{equation}
For each fixed background $\varphi$, this is expectation over
$\psi=\varphi+\zeta$ with density proportional to
$e^{-F(\psi)-\frac12(\psi-\varphi,\Sigma^{-1}(\psi-\varphi))}$.
We call this the conditional law at $\varphi$.
Write
$\operatorname{Cov}_{F,\Sigma}(f,g)
=T_{F,\Sigma}(fg)-(T_{F,\Sigma}f)(T_{F,\Sigma}g)$.
These are the Polchinski operators of~\cite[Definition~3.2 and Proposition~3.5]{BBD2024}.
For our scale step,
$\mathbf E_\lambda h=T_{N+\lambda R,C_T}h(Q_T\varphi')$.
If $F_r(\varphi)=-\log\mathbf E_{r\Sigma}e^{-F(\varphi+\zeta)}$,
Gaussian covariance addition gives
$T_{F,\Sigma}=T_{F_r,(1-r)\Sigma}T_{F,r\Sigma}$.

The background $\varphi$ occurs only in the Gaussian density when
we integrate in $\psi$. We can therefore differentiate a conditional
average without differentiating its observable. For bounded $h$,
differentiating the numerator and denominator gives
\[
 \begin{aligned}
 \partial_wT_{F,\Sigma}h
 &=T_{F,\Sigma}[h(\Sigma^{-1}w,\psi-\varphi)]
    -(T_{F,\Sigma}h)T_{F,\Sigma}(\Sigma^{-1}w,\psi-\varphi)\\
 &=\operatorname{Cov}_{F,\Sigma}(h,(\Sigma^{-1}w,\psi)).
 \end{aligned}
\]
If $\Hess F\ge-\Sigma^{-1}/2$, the Hessian of the conditional
Hamiltonian is at least $\Sigma^{-1}/2$. Brascamp--Lieb~\cite{BL1976}
then gives
\[
 \operatorname{Var}_{F,\Sigma}((\Sigma^{-1}w,\psi))
 \le T_{F,\Sigma}(\Sigma^{-1}w,
              [\Sigma^{-1}+\Hess F]^{-1}\Sigma^{-1}w)
 \le2(w,\Sigma^{-1}w).
\]
Using $\operatorname{Var}_{F,\Sigma}(h)\le\|h\|_\infty^2$
in Cauchy--Schwarz proves
\begin{equation}\label{cr:conditional-derivative}
 |\partial_wT_{F,\Sigma}h|
 \le\sqrt2\|h\|_\infty\sqrt{(w,\Sigma^{-1}w)}.
\end{equation}

The next lemma lets us differentiate a covariance using only first
derivatives of its observables. This is what we need when the
observables are a local term and a derivative of the interaction.
The proof splits the Gaussian integration so that the background
derivative falls on the outer density. The convolved-Hessian estimate
keeps the corresponding conditional Hamiltonian convex.
Related weighted covariance estimates appear
in~\cite[Theorem~3.1 and Remark~3.2]{Menz2014}.
\begin{lemma}\label{lem:conditional-covariance}
Let $F\in C^2(\R^n)$, $\Sigma>0$, and let $W_1,W_2$ be positive
diagonal matrices such that
\[
 \|W_i\Sigma W_i^{-1}\|_{\mathrm{op}}
 \|W_i\Hess F W_i^{-1}\|_{\infty,\mathrm{op}}\le\tfrac18,
 \qquad i=1,2.
\]
For $C^1$ functions $f,g$ with the indicated bounded gradients,
every direction $w$, $0\le r\le1$ and $i=1,2$, uniformly in the
background,
\begin{align}
 |\operatorname{Cov}_{F,\Sigma}(f,g)|
 &\le2\|W_1^{-1}\Sigma W_2^{-1}\|_{\mathrm{op}}
          \|W_1\nabla f\|_{\infty,2}\|W_2\nabla g\|_{\infty,2},
          \label{cr:product-covariance}\\
 \|W_i\Hess F_rW_i^{-1}\|_{\infty,\mathrm{op}}
 &\le2\|W_i\Hess F W_i^{-1}\|_{\infty,\mathrm{op}},
          \label{eq:hessian-convolution}\\
 |\partial_w\operatorname{Cov}_{F,\Sigma}(f,g)|
 &\le8\sqrt2\|W_1^{-1}\Sigma W_2^{-1}\|_{\mathrm{op}}
          \sqrt{(w,\Sigma^{-1}w)}
          \|W_1\nabla f\|_{\infty,2}\|W_2\nabla g\|_{\infty,2}.
          \label{cr:covariance-derivative}
\end{align}
\end{lemma}
\begin{proof}
Fix the background and write $\mathbf E$ for conditional expectation.
We first check convexity of the conditional Hamiltonian. For the
symmetric matrices $\Sigma$ and $\Hess F(\psi)$, the operator norm
is the largest modulus of an eigenvalue. Conjugation by $W_i$
preserves these eigenvalues, so the conjugated norms bound the
original ones. The hypothesis therefore gives
\[
 \|\Sigma\|_{\mathrm{op}}\|\Hess F\|_{\infty,\mathrm{op}}
 \le\|W_i\Sigma W_i^{-1}\|_{\mathrm{op}}
       \|W_i\Hess F W_i^{-1}\|_{\infty,\mathrm{op}}
 \le\tfrac18.
\]
In particular, $\id\le\|\Sigma\|_{\mathrm{op}}\Sigma^{-1}$ and
\[
 \begin{aligned}
 \Sigma^{-1}+\Hess F
 &\ge\Sigma^{-1}-\|\Hess F\|_{\infty,\mathrm{op}}\id\\
 &\ge\bigl(1-\|\Sigma\|_{\mathrm{op}}
                  \|\Hess F\|_{\infty,\mathrm{op}}\bigr)\Sigma^{-1}
 \ge\tfrac78\Sigma^{-1}.
 \end{aligned}
\]

Consider the operator
\[
 \mathcal L=\Tr(\Sigma\nabla^2)
            -(\psi-\varphi+\Sigma\nabla F,\nabla).
\]
Integration by parts against the conditional density gives, for
smooth $h,k$,
\[
 \begin{aligned}
 \mathbf E[(\mathcal Lh)k]
 &=\sum_{i,j}\Sigma_{ij}\mathbf E\left[
       k\partial_i\partial_jh
       -k\bigl(\partial_iF+[\Sigma^{-1}(\psi-\varphi)]_i\bigr)
           \partial_jh\right]\\
 &=-\sum_{i,j}\Sigma_{ij}\mathbf E[(\partial_i k)(\partial_j h)]
 =-\mathbf E(\nabla h,\Sigma\nabla k).
 \end{aligned}
\]
The curvature bound gives a Poincar\'e inequality by
Brascamp--Lieb. Hence the mean-zero weak solution $u$ of
$-\mathcal Lu=f-\mathbf E f$ exists, and the preceding identity gives
the covariance representation~\cite[Section~2, (2.5)--(2.9)]{Menz2014},
\[
 \operatorname{Cov}_{F,\Sigma}(f,g)
 =\mathbf E[(-\mathcal Lu)g]
 =\mathbf E(\nabla u,\Sigma\nabla g).
\]
The covariance is therefore controlled by a weighted $L^2$ bound
on $\nabla u$.
Differentiating the Poisson equation once gives
\[
 (\id-\mathcal L+\Hess F\Sigma)\nabla u=\nabla f.
\]
For either weight $W=W_i$, the hypothesis implies
\[
 \|W\Hess F\Sigma W^{-1}\|_{\infty,\mathrm{op}}
 \le\|W\Hess F W^{-1}\|_{\infty,\mathrm{op}}
       \|W\Sigma W^{-1}\|_{\mathrm{op}}\le\tfrac18.
\]
Test the gradient equation against $W^2\nabla u$ and integrate.
The diffusion contributes a nonnegative term, while the displayed
bound absorbs the Hessian term,
\[
 \begin{aligned}
 \mathbf E(W\nabla u,W\nabla f)
 &=\mathbf E\|W\nabla u\|_2^2
   +\sum_k\mathbf E(\nabla(W\nabla u)_k,
                         \Sigma\nabla(W\nabla u)_k)\\
 &{}+\mathbf E(W\nabla u,W\Hess F\Sigma W^{-1}W\nabla u)\\
 &\ge\tfrac78\mathbf E\|W\nabla u\|_2^2.
 \end{aligned}
\]
For the weak solution this test is justified by local elliptic
regularity and a smooth cutoff in field space. The cutoff term tends
to zero since $\mathbf E\|\nabla u\|_2^2<\infty$.
Cauchy--Schwarz now gives
$(\mathbf E\|W\nabla u\|_2^2)^{1/2}\le2\|W\nabla f\|_{\infty,2}$.
Putting this in the covariance representation proves
\[
 \begin{aligned}
 |\operatorname{Cov}_{F,\Sigma}(f,g)|
 &=|\mathbf E(W_1\nabla u,
             W_1^{-1}\Sigma W_2^{-1}W_2\nabla g)|\\
 &\le\|W_1^{-1}\Sigma W_2^{-1}\|_{\mathrm{op}}
       (\mathbf E\|W_1\nabla u\|_2^2)^{1/2}
       \|W_2\nabla g\|_{\infty,2}\\
 &\le2\|W_1^{-1}\Sigma W_2^{-1}\|_{\mathrm{op}}
       \|W_1\nabla f\|_{\infty,2}\|W_2\nabla g\|_{\infty,2}.
 \end{aligned}
\]

To bound the convolved Hessian, differentiate with respect to the
background while keeping $\zeta$ as the integration variable.
The quotient rule gives
\[
 \nabla T_{F,r\Sigma}f
 =T_{F,r\Sigma}\nabla f
       -\operatorname{Cov}_{F,r\Sigma}(f,\nabla F).
\]
We apply the covariance bound with reciprocal weights. If $W$ is
one of the weights in the hypothesis, symmetry and transposition give
\[
 \|W^{-1}\Sigma W\|_{\mathrm{op}}
 =\|W\Sigma W^{-1}\|_{\mathrm{op}},\qquad
 \|W^{-1}\Hess F W\|_{\infty,\mathrm{op}}
 =\|W\Hess F W^{-1}\|_{\infty,\mathrm{op}}.
\]
Thus $W^{-1}$ satisfies the same hypothesis. For $\|v\|_2=1$,
use the weights $W,W^{-1}$ and
$\nabla(Wv,\nabla F)=\Hess F W v$ to obtain
\[
 \begin{aligned}
 |(v,W\operatorname{Cov}_{F,r\Sigma}(f,\nabla F))|
 &=|\operatorname{Cov}_{F,r\Sigma}(f,(Wv,\nabla F))|\\
 &\le2r\|W^{-1}\Sigma W\|_{\mathrm{op}}
        \|W\nabla f\|_{\infty,2}
        \|W^{-1}\Hess F W v\|_{\infty,2}\\
 &\le2r\|W\Sigma W^{-1}\|_{\mathrm{op}}
        \|W\Hess F W^{-1}\|_{\infty,\mathrm{op}}
        \|W\nabla f\|_{\infty,2}.
 \end{aligned}
\]
Taking the supremum over $v$ and the background proves
\[
 \begin{aligned}
 \|W\nabla T_{F,r\Sigma}f\|_{\infty,2}
 &\le\bigl(1+2r\|W\Sigma W^{-1}\|_{\mathrm{op}}
                   \|W\Hess F W^{-1}\|_{\infty,\mathrm{op}}\bigr)
                    \|W\nabla f\|_{\infty,2}\\
 &\le2\|W\nabla f\|_{\infty,2}.
 \end{aligned}
\]
Since $\nabla F_r=T_{F,r\Sigma}\nabla F$, apply this estimate to
$(W^{-1}v,\nabla F)$. For every $\|v\|_2=1$ this gives
\[
 \begin{aligned}
 \|W\Hess F_rW^{-1}v\|_{\infty,2}
 &=\|W\nabla T_{F,r\Sigma}(W^{-1}v,\nabla F)\|_{\infty,2}\\
 &\le2\|W\Hess F W^{-1}v\|_{\infty,2}
 \le2\|W\Hess F W^{-1}\|_{\infty,\mathrm{op}},
 \end{aligned}
\]
which proves~\eqref{eq:hessian-convolution}.

It remains to differentiate the covariance. We split the Gaussian
integration, keeping the gradients of $f,g$ in the inner average
and differentiating the outer density with respect to the background. Put
$f_r=T_{F,r\Sigma}f$ and $g_r=T_{F,r\Sigma}g$.
The Gaussian-convolution proof of the evolution identities
in~\cite[(3.25)--(3.26)]{BBD2024} applies here because the conditional
Hamiltonians have the quadratic lower bound proved above. It gives
\[
 \partial_rf_r=\tfrac12\Tr(\Sigma\nabla^2f_r)
                     -(\nabla F_r,\Sigma\nabla f_r),
\]
and the same equation for $g_r$. When differentiating
$T_{F_r,(1-r)\Sigma}(f_rg_r)$, we also differentiate $F_r$ and the
remaining covariance $(1-r)\Sigma$. Their contribution cancels
the separate evolution terms of $f_r$ and $g_r$. The product rule
leaves the cross derivative,
\[
 \begin{aligned}
 \frac d{dr}T_{F_r,(1-r)\Sigma}(f_rg_r)
 &=T_{F_r,(1-r)\Sigma}\left[
       \partial_r(f_rg_r)-\tfrac12\Tr(\Sigma\nabla^2(f_rg_r))
       +(\nabla F_r,\Sigma\nabla(f_rg_r))\right]\\
 &=\tfrac12T_{F_r,(1-r)\Sigma}\left[
       g_r\Tr(\Sigma\nabla^2 f_r)+f_r\Tr(\Sigma\nabla^2 g_r)
       -\Tr(\Sigma\nabla^2(f_rg_r))\right]\\
 &=-T_{F_r,(1-r)\Sigma}(\nabla f_r,\Sigma\nabla g_r).
 \end{aligned}
\]
At $r=0$ the average is $T_{F,\Sigma}(fg)$; at $r=1$ it is
$(T_{F,\Sigma}f)(T_{F,\Sigma}g)$. Integration therefore gives
\begin{equation}\label{cr:covariance-splitting}
 \operatorname{Cov}_{F,\Sigma}(f,g)
 =\int_0^1T_{F_r,(1-r)\Sigma}(\nabla f_r,\Sigma\nabla g_r)\,dr.
\end{equation}
The conditional-gradient estimate bounds the observable by
\[
 \begin{aligned}
 \|(\nabla f_r,\Sigma\nabla g_r)\|_\infty
 &\le\|W_1^{-1}\Sigma W_2^{-1}\|_{\mathrm{op}}
        \|W_1\nabla f_r\|_{\infty,2}\|W_2\nabla g_r\|_{\infty,2}\\
 &\le4\|W_1^{-1}\Sigma W_2^{-1}\|_{\mathrm{op}}
        \|W_1\nabla f\|_{\infty,2}\|W_2\nabla g\|_{\infty,2}.
 \end{aligned}
\]
For the outer conditional law,
\eqref{eq:hessian-convolution} with unit weights gives
\[
 \Hess F_r\ge-2\|\Hess F\|_{\infty,\mathrm{op}}\id
 \ge-\tfrac14\Sigma^{-1},\qquad
 [(1-r)\Sigma]^{-1}+\Hess F_r
 \ge\tfrac12[(1-r)\Sigma]^{-1}.
\]
For each $r$, the function $(\nabla f_r,\Sigma\nabla g_r)$ depends
on the integration variable, not on the outer background. We can
therefore apply~\eqref{cr:conditional-derivative} to the outer density
without taking another derivative of $f_r$ or $g_r$, obtaining
\[
 \begin{aligned}
 |\partial_w\operatorname{Cov}_{F,\Sigma}(f,g)|
 &\le\sqrt2\sqrt{(w,\Sigma^{-1}w)}
        \int_0^1\frac{\|(\nabla f_r,\Sigma\nabla g_r)\|_\infty}
                           {\sqrt{1-r}}\,dr\\
 &\le8\sqrt2\|W_1^{-1}\Sigma W_2^{-1}\|_{\mathrm{op}}
        \sqrt{(w,\Sigma^{-1}w)}
        \|W_1\nabla f\|_{\infty,2}\|W_2\nabla g\|_{\infty,2}.
 \end{aligned}
\]
The integrable factor $(1-r)^{-1/2}$ justifies differentiation under
the integral. The endpoint formulas for the undifferentiated
interpolation follow by continuity. Mollification in field space
preserves the weighted Hessian and gradient bounds.
For sufficiently small mollification radii, the conditional
Hamiltonians have a common quadratic lower bound, so dominated
convergence passes the estimates to the stated $C^2,C^1$ inputs.
The constants in the estimates are unchanged.
\end{proof}

For the scale step, take $F=N+\lambda R$, $0\le\lambda\le1$,
and $\Sigma=C_T$. Let $e_i$ now denote the unit vector in the field
coordinate $i$, so $(Q_Te_i)_a=Q_T(a,i)$, and let $W_x$ be diagonal
with $(W_x)_{aa}=e^{\rho d(a,x)/\ell'}$.
Assume $0<\rho\le\rho_0$ and $\rho/\ell'\le\mu/(8\ell)$.
The triangle inequality gives
\[
 |(W_x\Hess F W_x^{-1})_{ab}|
 \le e^{\rho d(a,b)/\ell'}\sup_\varphi|\partial_a\partial_bF|.
\]
The expression on the right is symmetric in $a,b$. Its row sums,
and hence its column sums, are at most $H\zeta_\ell\ell^{-2}$ by
\eqref{ap:source-estimates}. Schur's inequality therefore gives
$\|W_x\Hess F W_x^{-1}\|_{\infty,\mathrm{op}}
\le H\zeta_\ell\ell^{-2}$.
The same argument for $C_T$, using~\eqref{cr:integrated-operators},
verifies the hypothesis of Lemma~\ref{lem:conditional-covariance},
\[
 \|W_xC_TW_x^{-1}\|_{\mathrm{op}}
 \|W_x\Hess F W_x^{-1}\|_{\infty,\mathrm{op}}
 \le C\ell'^2H\zeta_\ell\ell^{-2}\le\tfrac18.
\]
For observables localized near $x,y$, the covariance bound retains
their separation,
\[
 \begin{aligned}
 \|W_x^{-1}C_TW_y^{-1}\|_{\mathrm{op}}
 &\le\|W_x^{-1}C_TW_x\|_{\mathrm{op}}
       \sup_a e^{-\rho[d(a,x)+d(a,y)]/\ell'}\\
 &\le C_\rho\ell'^2e^{-\rho d(x,y)/\ell'}.
 \end{aligned}
\]
The background derivative in~\eqref{cr:covariance-derivative}
introduces an inverse covariance. Its square root contributes
$\ell'^{-2}$ even if only the last half of the Gaussian integration
remains. Indeed, for $0\le S\le T/2$,
\begin{equation}\label{eq:fresh-gaussian-cost}
 \begin{aligned}
 (Q_Te_j,(C_T-C_S)^{-1}Q_Te_j)
 &=\left[\frac A{e^{(T-S)A}-1}\right]_{jj}\\
 &\le\frac CT(e^{-TA/4})_{jj}\le C\ell'^{-4}.
 \end{aligned}
\end{equation}
The identity follows from
$C_T-C_S=A^{-1}e^{-SA}(\id-e^{-(T-S)A})$ and $Q_T^2=e^{-TA}$.
For the inequality use
$a/(e^{(T-S)a}-1)\le CT^{-1}e^{-Ta/4}$ for $a>0$, then the diagonal
heat bound. At $S=0$, the square root in
\eqref{cr:covariance-derivative} cancels the factor $\ell'^2$ from
the covariance. Once the smallness condition holds, these constants
are independent of $M$, cutoff and volume.

We shall also need the gradient of $(Q_Te_i,\nabla F)$.
The heat-kernel bound and its spatial sum give
\[
 \begin{aligned}
 \|W_iQ_Te_i\|_2^2
 &=\sum_a e^{2\rho d(a,i)/\ell'}Q_T(a,i)^2\\
 &\le\left(\sup_a e^{\rho d(a,i)/\ell'}Q_T(a,i)\right)
       \sum_a e^{\rho d(a,i)/\ell'}Q_T(a,i)
 \le C_\rho\ell'^{-2}.
 \end{aligned}
\]
Combining this with the weighted Hessian bound yields
\begin{equation}\label{eq:weighted-heat-column}
 \begin{aligned}
 \|W_i\Hess FQ_Te_i\|_{\infty,2}
 &\le\|W_i\Hess F W_i^{-1}\|_{\infty,\mathrm{op}}
         \|W_iQ_Te_i\|_2
 \le C_\rho H\zeta_\ell\ell^{-2}\ell'^{-1}.
 \end{aligned}
\end{equation}
Similarly, for $\rho\le6\mu$, squaring~\eqref{ap:ordinary-gradient}
gives
\[
 \|W_x\nabla g_x\|_{\infty,2}^2
 \le C(K+1)^2\zeta_\ell^2\ell^{-4}
       \sum_a e^{-6\mu d(a,x)/\ell}
 \le C(K+1)^2\zeta_\ell^2\ell^{-2}.
\]
These sums involve the diffusion area, rather than the number of sites.

The first derivative of $G_x$ contains a covariance with
$[Q_T\nabla F]_i$. The Hessian of $G_x$ contains the background derivative
of this covariance. Write
$\operatorname{Cov}_\lambda$ for covariance under $\mathbf E_\lambda$,
and take $\rho=6\mu$, which is admissible for $M\ge48$.
By~\eqref{cr:product-covariance} and the preceding gradient bounds,
\begin{equation}\label{fc:force-covariance}
 \begin{aligned}
 |\operatorname{Cov}_\lambda(g_x,[Q_T\nabla F]_i)|
 &\le C\ell'^2e^{-6\mu d(i,x)/\ell'}
       [(K+1)\zeta_\ell\ell^{-1}]
       [H\zeta_\ell\ell^{-2}\ell'^{-1}]\\
 &\le C_{M,K}H\zeta_\ell^2\ell'^{-2}e^{-6\mu d(i,x)/\ell'}.
 \end{aligned}
\end{equation}
A derivative in the outgoing background $\varphi'_j$ is a derivative
of the conditional law in direction $Q_Te_j$.
Thus~\eqref{cr:covariance-derivative} has the same product of weighted
gradients, with the additional factor bounded by $C\ell'^{-2}$ in
\eqref{eq:fresh-gaussian-cost}. We obtain
\begin{equation}\label{fc:force-covariance-derivative}
 |\partial_j\operatorname{Cov}_\lambda(g_x,[Q_T\nabla F]_i)|
 \le C_{M,K}H\zeta_\ell^2\ell'^{-4}e^{-6\mu d(i,x)/\ell'}.
\end{equation}
Here $F=N+\lambda R$, and both bounds hold uniformly in
$\lambda\in[0,1]$ and the background, with $C_{M,K}$ independent of
$H,B$, cutoff and volume.

\subsection{The neutral convolution}
The Hessian of the neutral convolution has an averaged Hessian term
and a covariance term. The latter is quadratic in $H\zeta_\ell$.
For the former, zero row sums supply two spatial differences, giving
a linear bound with coefficient independent of $M$.

\begin{proposition}\label{prop:neutral-smoothing}
Suppose $N$ is invariant under constant shifts and satisfies the
row bound~\eqref{ap:full-hessian-row} with $N$ in place of $N+R$.
For $M\ge16$, in the heat regime and under the standing smallness
conditions, its Gaussian convolution is neutral and
\begin{equation}\label{ap:neutral-recursion}
 \begin{aligned}
 &\sup_i\sum_j e^{\mu d(i,j)/(2\ell')}
   \sup_{\varphi'}\left|\partial_i\partial_j
      [-\log\mathbf E_{C_T}e^{-N(Q_T\varphi'+\zeta)}]\right|\\
 &\le[C_\mu H+C_{\beta,\mu,M}H^2\zeta_\ell]
          \zeta_\ell\ell'^{-2}.
 \end{aligned}
\end{equation}
The constant $C_\mu$ is independent of $M,K,H,B$.
\end{proposition}
\begin{proof}
Differentiate $N(\varphi+s\mathbf1)=N(\varphi)$ first in $s$
and then in $\varphi_a$. This gives
$\sum_b\partial_a\partial_bN=0$. The column sums also vanish,
since $\partial_a\partial_bN=\partial_b\partial_aN$. In the following expansion, the terms containing
$u_aw_a$ and $u_bw_b$ therefore sum to zero,
\[
 \begin{aligned}
 \sum_{a,b}\partial_a\partial_bN(u_a-u_b)(w_a-w_b)
 &=-\sum_{a,b}\partial_a\partial_bN(u_aw_b+u_bw_a)\\
 &=-2(u,\Hess Nw).
 \end{aligned}
\]
Equivalently,
\begin{equation}\label{ap:neutral-difference-form}
 (u,\Hess Nw)=-\tfrac12\sum_{a,b}\partial_a\partial_bN
                                     (u_a-u_b)(w_a-w_b).
\end{equation}
The conditional Hessian identity~\cite[Lemma~3.12]{BBD2024} is
\[
 \partial_i\partial_j[-\log\mathbf E_{C_T}e^{-N(Q_T\varphi'+\zeta)}]
 =T_{N,C_T}(Q_Te_i,\Hess NQ_Te_j)
       -\operatorname{Cov}_{N,C_T}
                    ((Q_Te_i,\nabla N),(Q_Te_j,\nabla N)).
\]
The averages have background $Q_T\varphi'$. For the first term,
the heat-difference bound gives
\[
 |Q_T(a,i)-Q_T(b,i)|
 \le Cd(a,b)\ell'^{-3}e^{-2\mu d(a,i)/\ell'+2\mu d(a,b)/\ell'}.
\]
For $M\ge16$ the second spatial moment costs only $C_\mu\ell^2$,
independently of $M$, since
\[
 d(a,b)^2e^{4\mu d(a,b)/\ell'}
 \le\bigl[d(a,b)^2e^{-\mu d(a,b)/(4\ell)}\bigr]
           e^{\mu d(a,b)/(2\ell)}
 \le C_\mu\ell^2 e^{\mu d(a,b)/(2\ell)}.
\]
Thus~\eqref{ap:neutral-difference-form} and the weighted row bound give
\[
 \begin{aligned}
 |(Q_Te_i,\Hess NQ_Te_j)|
 &\le C\ell'^{-6}\sum_a e^{-2\mu[d(a,i)+d(a,j)]/\ell'}
       \sum_b d(a,b)^2e^{4\mu d(a,b)/\ell'}|\partial_a\partial_bN|\\
 &\le C_\mu\ell^2\ell'^{-6}\sum_a e^{-2\mu[d(a,i)+d(a,j)]/\ell'}
       \sum_b e^{\mu d(a,b)/(2\ell)}\sup_\varphi|\partial_a\partial_bN|\\
 &\le C_\mu H\zeta_\ell\ell'^{-6}
                     \sum_a e^{-2\mu[d(a,i)+d(a,j)]/\ell'}\\
 &\le C_\mu H\zeta_\ell\ell'^{-4}e^{-\mu d(i,j)/\ell'}.
 \end{aligned}
\]
The last sum is bounded by
\[
 \sum_a e^{-2\mu[d(a,i)+d(a,j)]/\ell'}
 \le e^{-\mu d(i,j)/\ell'}\sum_a e^{-\mu d(a,i)/\ell'}
 \le C_\mu\ell'^2e^{-\mu d(i,j)/\ell'}.
\]
For the covariance use $(W_i)_{aa}=e^{\mu d(a,i)/\ell'}$.
The two observables have gradients $\Hess NQ_Te_i$ and
$\Hess NQ_Te_j$. Thus~\eqref{cr:product-covariance}
and~\eqref{eq:weighted-heat-column} with $F=N$ give
\[
 \begin{aligned}
 |\operatorname{Cov}_{N,C_T}((Q_Te_i,\nabla N),(Q_Te_j,\nabla N))|
 &\le C\ell'^2e^{-\mu d(i,j)/\ell'}
                (H\zeta_\ell\ell^{-2}\ell'^{-1})^2\\
 &\le C_{\beta,\mu,M}H^2\zeta_\ell^2\ell'^{-4}e^{-\mu d(i,j)/\ell'}.
 \end{aligned}
\]
Multiply these two bounds by $e^{\mu d(i,j)/(2\ell')}$ and sum
over $j$. Since $\sum_j e^{-\mu d(i,j)/(2\ell')}\le C_\mu\ell'^2$,
we obtain~\eqref{ap:neutral-recursion}. Neutrality is preserved by
Gaussian convolution and by composition with $Q_T$, which maps
constants to constants.
\end{proof}

\subsection{Gaussian averaging}
We compare Gaussian averaging with an independent heat average of
the phase. The error involves only derivatives with the phase held
fixed and vanishes on the original cosine.

\begin{lemma}\label{lem:gaussian-phase-damping}
Let $f$ be $C^1$, periodic under constant shifts with period $2\pi/c_t$
and with zero period average. Suppose the norms below are finite. Fix $x\in\Lambda$ and $M\ge2$.
In the heat regime,
\begin{equation}\label{ap:periodic-value-general}
 \|\mathbf E_{C_T}f(\cdot+\zeta)\|_*
 \le e^{-\frac{\beta}2(C_{t+T}-C_t)_{xx}}
 \left[\|f\|_*+C_{\beta,\mu}\ell^2
        \sup_i e^{\mu d(i,x)/\ell}\|\partial_i^{(x)}f\|_*\right].
\end{equation}
The constant is independent of $M$, cutoff and volume.
\end{lemma}
\begin{proof}
Split the covariance as $C_T=C_s+(C_T-C_s)$, and let
$\zeta,\eta$ be independent centred Gaussian vectors with these
covariances. We leave $\zeta$ unchanged. From $\eta$ we subtract a
constant field which removes its contribution to the phase, and
replace that contribution by an independent heat average in the
angle. Varying $s$ from zero to $T$ then restores the original
correlation between the phase and the other field variables. Thus put
\[
 f_s(\varphi)=
 e^{\frac{\beta}2(C_{t+T}-C_{t+s})_{xx}\Theta_t^2}
 \mathbf E f\bigl(\varphi+\zeta+\eta
          -\sqrt\beta(Q_t\eta)_xc_t^{-1}\mathbf1\bigr).
\]
The subtraction removes the change in the phase at $x$,
\[
 \sqrt\beta\bigl[Q_t(\eta
       -\sqrt\beta(Q_t\eta)_xc_t^{-1}\mathbf1)\bigr]_x
 =\sqrt\beta(Q_t\eta)_x
       -\sqrt\beta(Q_t\eta)_xc_t^{-1}c_t=0.
\]
The heat semigroup on the circle restores the phase variance,
$\beta[Q_t(C_T-C_s)Q_t]_{xx}
=\beta(C_{t+T}-C_{t+s})_{xx}$, independently of $\eta$.
At the two endpoints we therefore have
\[
 f_T=\mathbf E_{C_T}f(\cdot+\zeta),\qquad
 \|f_0\|_*\le
 e^{-\frac{\beta}2(C_{t+T}-C_t)_{xx}}\|f\|_*.
\]

To estimate $f_T-f_0$, differentiate the commuting averages.
As $s$ increases, the covariance $C_s$ grows at rate $e^{-sA}$.
The covariance of $\eta$ decreases at the same rate, with derivatives
$\partial_i^{(x)}$ in place of $\partial_i$. Differentiating also
the remaining angular variance gives
\[
 \partial_s f_s=
 \left[\frac12\sum_{i,j}(e^{-sA})_{ij}
       (\partial_i\partial_j-\partial_i^{(x)}\partial_j^{(x)})
       -\frac{\beta}2(e^{-(t+s)A})_{xx}\Theta_t^2\right]f_s.
\]
The terms with two angular derivatives cancel. Indeed,
\[
 \begin{aligned}
 \partial_i\partial_j-\partial_i^{(x)}\partial_j^{(x)}
 &=\sqrt\beta\Theta_t
    [Q_t(i,x)\partial_j^{(x)}+Q_t(j,x)\partial_i^{(x)}]
   +\beta Q_t(i,x)Q_t(j,x)\Theta_t^2,\\
 \sum_{i,j}(e^{-sA})_{ij}Q_t(i,x)Q_t(j,x)
 &=(Q_te^{-sA}Q_t)_{xx}=(e^{-(s+t)A})_{xx},\\
 \partial_s f_s
 &=\sqrt\beta\Theta_t\sum_i(e^{-(s+t/2)A})_{ix}
                         \partial_i^{(x)}f_s.
 \end{aligned}
\]
The phase-fixed derivative is controlled by the hypothesis. To bound
the angular derivative without differentiating $f$ again, we take it
on the remaining heat kernel on the circle. The derivative of this kernel
has $L^1$ norm at most $Ca^{-1/2}$
for $a\le1$, while for $a\ge1$ the Fourier bound
$\sum_{q\ne0}|q|e^{-q^2a}\le Ce^{-a}$ applies. Commuting with the
semigroup in~\eqref{ap:circle-norm}, we obtain
\[
 \|\Theta_t e^{a\Theta_t^2}h\|_*
 \le C(1+a^{-1/2})e^{-a}\|h\|_*.
\]
Apply this bound to the angular heat semigroup in $f_s$ and contract
the other averages. Integrating from the independent-phase endpoint
to the original average gives
\begin{equation}\label{all:duhamel}
 \begin{aligned}
 \|f_T-f_0\|_*
 &\le C\int_0^T
 \left(1+\frac1{\sqrt{\beta(C_{t+T}-C_{t+s})_{xx}}}\right)
 e^{-\frac{\beta}2(C_{t+T}-C_{t+s})_{xx}}\\
 &{}\cdot
 \left\|\sum_i\sqrt\beta(e^{-(s+t/2)A})_{ix}
                         \partial_i^{(x)}f\right\|_*\,ds.
 \end{aligned}
\end{equation}

We can improve the time decay in~\eqref{all:duhamel} because the
phase-fixed derivatives sum to zero,
\[
 \sum_i\partial_i^{(x)}f
 =\sum_i\partial_if-\sqrt\beta\sum_iQ_t(i,x)\Theta_tf
 =c_t\Theta_tf-c_t\Theta_tf=0.
\]
Subtracting the same value $(e^{-(s+t/2)A})_{xx}$ from every
coefficient therefore leaves the sum unchanged. By evenness of the
kernel about $x$, the resulting difference is half a centred second
spatial difference. Telescoping twice along a shortest path and using
Lemma~\ref{cr:heat-bounds} gives
\[
 \begin{aligned}
 |(e^{-(s+t/2)A})_{ix}-(e^{-(s+t/2)A})_{xx}|
 &=\tfrac12|(e^{-(s+t/2)A})_{ix}
       +(e^{-(s+t/2)A})_{2x-i,x}-2(e^{-(s+t/2)A})_{xx}|\\
 &\le C\frac{d(i,x)^2}{(\ell^2+s)^2},\\
 \left\|\sum_i\sqrt\beta(e^{-(s+t/2)A})_{ix}
                         \partial_i^{(x)}f\right\|_*
 &=\left\|\sum_i\sqrt\beta
       [(e^{-(s+t/2)A})_{ix}-(e^{-(s+t/2)A})_{xx}]
                         \partial_i^{(x)}f\right\|_*\\
 &\le\frac{C_\beta}{(\ell^2+s)^2}
       \sum_i d(i,x)^2\|\partial_i^{(x)}f\|_*\\
 &\le\frac{C_{\beta,\mu}\ell^4}{(\ell^2+s)^2}
       \sup_i e^{\mu d(i,x)/\ell}\|\partial_i^{(x)}f\|_*.
 \end{aligned}
\]
Here
$\sum_i d(i,x)^2e^{-\mu d(i,x)/\ell}
\le4\sum_{n\ge0}(n+1)n^2e^{-\mu n/\ell}\le C_\mu\ell^4$.
The factor $(\ell^2+s)^{-2}$ must compensate for the phase damping
lost up to time $s$. The covariance logarithm bounds that loss, and
the lower heat diagonal controls the remaining circle variance,
\[
 \begin{aligned}
 e^{\frac{\beta}2(C_{t+s}-C_t)_{xx}}
 &\le C_\beta(1+s/\ell^2)^{\beta/(8\pi)},\\
 \beta(C_{t+T}-C_{t+s})_{xx}
 &=\beta\int_{t+s}^{t+T}(e^{-uA})_{xx}\,du
 \ge c_\beta\frac{T-s}{\ell'^2}.
 \end{aligned}
\]
Take $e^{-\beta(C_{t+T}-C_t)_{xx}/2}$ out of~\eqref{all:duhamel}
and set $u=1+s/\ell^2$, so $ds=\ell^2du$ and
$T-s=\ell^2(M^2-u)$. Then
\[
 \begin{aligned}
 &\int_0^T
 \left(1+\frac1{\sqrt{\beta(C_{t+T}-C_{t+s})_{xx}}}\right)
 e^{\frac{\beta}2(C_{t+s}-C_t)_{xx}}
       \frac{\ell^4}{(\ell^2+s)^2}\,ds\\
 &\le C_\beta\ell^2\int_1^{M^2}
        u^{\beta/(8\pi)-2}
        \left(1+\sqrt{\frac{M^2}{M^2-u}}\right)du.
 \end{aligned}
\]
The second spatial difference gives the power $\beta/(8\pi)-2$,
which is integrable at infinity. Near the upper endpoint $u=M^2$
the angular variance tends to zero, and differentiation of its heat kernel
costs the inverse square root of that variance. This singularity is
also integrable. To separate the two estimates, split at $M^2/2$
and set $v=u/M^2$ in the upper half,
\begin{equation}\label{ap:critical-integral}
 \begin{aligned}
 &\int_1^{M^2}u^{\beta/(8\pi)-2}
       \left(1+\sqrt{\frac{M^2}{M^2-u}}\right)du\\
 &\le(1+\sqrt2)\int_1^\infty u^{\beta/(8\pi)-2}\,du
    +C_\beta M^{\beta/(4\pi)-2}
       \int_{1/2}^1(1+(1-v)^{-1/2})\,dv\\
 &\le\frac{1+\sqrt2}{1-\beta/(8\pi)}
                 +C_\beta M^{\beta/(4\pi)-2}
 \le\frac{C_\beta}{1-\beta/(8\pi)}.
 \end{aligned}
\end{equation}
Substitution in~\eqref{all:duhamel} yields
\begin{equation}\label{ap:gaussian-phase-error}
 \|f_T-f_0\|_*
 \le C_{\beta,\mu}\ell^2e^{-\frac\beta2(C_{t+T}-C_t)_{xx}}
       \sup_i e^{\mu d(i,x)/\ell}\|\partial_i^{(x)}f\|_*.
\end{equation}
The bound on $f_0$ then proves~\eqref{ap:periodic-value-general}.
Convolution approximation gives the same conclusion for $C^1$
functions, since the bounds above use only the first derivatives
in the statement.
\end{proof}

To pass from the Gaussian average to $\mathbf E_0g_x$, we compare
at the actual conditional mean. A small Hessian controls fluctuations,
but not the displacement of this mean. The displacement is harmless
because the preceding estimate is uniform in the background.

\begin{lemma}\label{lem:conditional-gaussian-comparison}
In the heat regime, let $F\in C^2$,
$0<s\le T$ and $s\|\Hess F\|_{\infty,\mathrm{op}}\le1/2$.
For bounded $C^1$ functions $f$ with bounded first derivatives,
\begin{equation}\label{cr:centered-gaussian-comparison}
 \|T_{F,C_s}f\|_\infty
 \le\|\mathbf E_{C_s}f(\cdot+\zeta)\|_\infty
       +Cs\|\Hess F\|_{\infty,\mathrm{op}}\sum_i\|\partial_i f\|_\infty.
\end{equation}
At each background, let $\mathbf E$ denote its conditional expectation.
For each direction $w$, comparison at its actual mean gives
\begin{equation}\label{cr:centered-gaussian-first}
 |\mathbf E\partial_wf(\psi)
       -\mathbf E_{C_s}\partial_wf(\mathbf E\psi+\zeta)| \le Cs\|\Hess F\|_{\infty,\mathrm{op}}\sqrt{(w,C_s^{-1}w)}
                                      \sum_i\|\partial_i f\|_\infty.
\end{equation}
If $F$ is neutral and $f$ has period $2\pi/c_t$ under common shifts
and zero shift average, then
\[
 \|T_{F,C_s}f\|_*
 \le\|\mathbf E_{C_s}f(\cdot+\zeta)\|_*
       +Cs\|\Hess F\|_{\infty,\mathrm{op}}\sum_i\|\partial_i f\|_*.
\]
The constants are independent of $M$, cutoff and volume.
\end{lemma}
\begin{proof}
Fix the background $\varphi$. Since $C_s\le s\id$ and
$s\|\Hess F\|_{\infty,\mathrm{op}}\le1/2$, the conditional Hamiltonian has
Hessian bounded below by
\[
 C_s^{-1}+\Hess F
 \ge C_s^{-1}-\|\Hess F\|_{\infty,\mathrm{op}}\id
 \ge(1-s\|\Hess F\|_{\infty,\mathrm{op}})C_s^{-1}
 \ge\tfrac12C_s^{-1}.
\]
Brascamp--Lieb gives
$\operatorname{Var}(h)\le2\mathbf E(\nabla h,C_s\nabla h)$.
Integrating the gradient of the conditional density gives the mean,
\[
 0=\mathbf E[\nabla F+C_s^{-1}(\psi-\varphi)],\qquad
 \mathbf E\psi-\varphi=-C_s\mathbf E\nabla F.
\]

Use the Ornstein--Uhlenbeck interpolation about this mean,
\cite[(2.7.3), (2.7.5)]{BGL2014},
\[
 f_u(\psi)=\mathbf E_{C_s}f\bigl(\mathbf E\psi+e^{-u}(\psi-\mathbf E\psi)
                            +\sqrt{1-e^{-2u}}\,\zeta\bigr).
\]
It joins $f$ to the constant $\mathbf E_{C_s}f(\mathbf E\psi+\zeta)$,
and $\|\partial_i f_u\|_\infty\le e^{-u}\|\partial_i f\|_\infty$.
Conditional integration by parts replaces the second derivatives of
$f_u$ by the gradient of the conditional Hamiltonian,
\[
 \mathbf E\Tr(C_s\Hess f_u)
 =\mathbf E(C_s\nabla F+\psi-\varphi,\nabla f_u).
\]
Substitution into the Ornstein--Uhlenbeck equation therefore gives
\[
 \begin{aligned}
 \frac d{du}\mathbf E f_u
 &=\mathbf E[\Tr(C_s\Hess f_u)-(\psi-\mathbf E\psi,\nabla f_u)]\\
 &=\mathbf E(C_s\nabla F+\mathbf E\psi-\varphi,\nabla f_u)\\
 &=\sum_i\operatorname{Cov}((C_se_i,\nabla F),\partial_i f_u).
 \end{aligned}
\]
Centering at the conditional mean has canceled the mean force.
Brascamp--Lieb now bounds its fluctuation by $\Hess F$. The image
of the coordinate vector $e_i$ under $C_s$ has a bound independent
of the number of sites,
\[
 \|C_se_i\|_2
 =\left\|\int_0^s e^{-vA}e_i\,dv\right\|_2
 \le\int_0^s\sqrt{(e^{-2vA})_{ii}}\,dv
 \le C\int_0^s(1+v)^{-1/2}\,dv\le C\sqrt s.
\]
Consequently,
\[
 \begin{aligned}
 \operatorname{Var}((C_se_i,\nabla F))
 &\le2\mathbf E(\Hess F C_se_i,C_s\Hess F C_se_i)\\
 &\le2s\|\Hess F\|_{\infty,\mathrm{op}}^2\|C_se_i\|_2^2\\
 &\le Cs^2\|\Hess F\|_{\infty,\mathrm{op}}^2.
 \end{aligned}
\]
Integrating the covariance identity and using Cauchy--Schwarz yields
\[
 \begin{aligned}
 |\mathbf E f-\mathbf E_{C_s}f(\mathbf E\psi+\zeta)|
 &\le\int_0^\infty\sum_i
       \sqrt{\operatorname{Var}((C_se_i,\nabla F))}
                   e^{-u}\|\partial_i f\|_\infty\,du\\
 &\le Cs\|\Hess F\|_{\infty,\mathrm{op}}\sum_i\|\partial_i f\|_\infty.
 \end{aligned}
\]
Taking the supremum over the background proves the scalar comparison.
To apply this calculation to $\partial_wf$, integrate the derivative
in direction $w$ against the Gaussian density. The observable then
still requires only one derivative,
\[
 \begin{aligned}
 \partial_i(\partial_wf)_u(\psi)
 &=\frac{e^{-u}}{\sqrt{1-e^{-2u}}}\mathbf E_{C_s}\bigl[
       (C_s^{-1}w,\zeta)\,
       \partial_i f(\mathbf E\psi+e^{-u}(\psi-\mathbf E\psi)
                          +\sqrt{1-e^{-2u}}\,\zeta)\bigr],\\
 \|\partial_i(\partial_wf)_u\|_\infty
 &\le\frac{e^{-u}}{\sqrt{1-e^{-2u}}}
       \sqrt{\mathbf E_{C_s}(C_s^{-1}w,\zeta)^2}\,
                   \|\partial_i f\|_\infty\\
 &=\frac{e^{-u}}{\sqrt{1-e^{-2u}}}
       \sqrt{(w,C_s^{-1}w)}\|\partial_i f\|_\infty.
 \end{aligned}
\]
Since $e^{-u}/\sqrt{1-e^{-2u}}$ is integrable over $(0,\infty)$,
the variance estimate proves~\eqref{cr:centered-gaussian-first}.
The same first-derivative bounds justify convolution approximation
for $C^1$ functions.
If $F$ is neutral, a change of variables in its conditional density gives
\[
 T_{F,C_s}[f(\cdot+\theta\mathbf1)](\varphi)
 =T_{F,C_s}f(\varphi+\theta\mathbf1).
\]
Integrating this identity against the heat kernel on the circle shows that
$T_{F,C_s}$ commutes with $e^{a\Theta_t^2}$. Apply the scalar estimate
to $e^{a\Theta_t^2}f$, multiply by $e^a$ and take the supremum over
$a\ge0$. Gaussian convolution and the first derivatives commute
with the same semigroup, giving the estimate in $\|\cdot\|_*$.
\end{proof}

\subsection{Estimates at the next scale}
\label{sec:coupling-propagation}\label{fc:section}
We first estimate the value by comparing $G_x$ with the average
$M^2\mathbf E_0g_x$ under $N$.
The effect of adding $\lambda R$ is quadratic in the coupling,
as differentiation of the numerator and normalization shows,
\[
 \frac d{d\lambda}\mathbf E_\lambda g_x
 =-\mathbf E_\lambda(g_xR)
       +(\mathbf E_\lambda g_x)(\mathbf E_\lambda R)
 =-\ell^{-2}\sum_y\operatorname{Cov}_\lambda(g_x,g_y).
\]
Integrating in $\lambda$ gives
\begin{equation}\label{fc:quadratic-identity}
 \begin{aligned}
 G_x-M^2\mathbf E_0g_x
 &=M^2\int_0^1(\mathbf E_\lambda g_x-\mathbf E_0g_x)\,d\lambda\\
 &=-M^2\ell^{-2}\sum_y\int_0^1(1-\lambda)
                     \operatorname{Cov}_\lambda(g_x,g_y)\,d\lambda.
 \end{aligned}
\end{equation}
The product covariance bound and local gradients give
\[
 \begin{aligned}
 |\operatorname{Cov}_\lambda(g_x,g_y)|
 &\le C_M\ell'^2\zeta_\ell^2\ell^{-2}e^{-6\mu d(x,y)/\ell'}
 \le C_M\zeta_\ell^2e^{-6\mu d(x,y)/\ell'},\\
 \|G_x-M^2\mathbf E_0g_x\|_\infty
 &\le C_MM^2\ell^{-2}\zeta_\ell^2\sum_y e^{-6\mu d(x,y)/\ell'}
 \le C_M\zeta_\ell^2.
 \end{aligned}
\]
For the average under $N$, the two errors in
Lemmas~\ref{lem:gaussian-phase-damping} and~\ref{lem:conditional-gaussian-comparison}
are, by~\eqref{ap:norms} and~\eqref{ap:ordinary-gradient},
\[
 \begin{aligned}
 \ell^2\sup_i e^{\mu d(i,x)/\ell}\|\partial_i^{(x)}g_x\|_*
 &\le CB\zeta_\ell^2,\\
 T\|\Hess N\|_{\infty,\mathrm{op}}\sum_i\|\partial_i g_x\|_*
 &\le C_{M,K}TH\zeta_\ell\ell^{-2}\zeta_\ell
 \le C_M\zeta_\ell^2.
 \end{aligned}
\]
Neutrality of $N$ makes $T_{N,C_T}$ commute with shift averaging.
The relation $Q_T\mathbf1=(c_{t+T}/c_t)\mathbf1$ identifies the
old and new periods, and invertibility of $Q_T$ preserves the field
supremum in $\|\cdot\|_*$. Thus $\Pi\mathbf E_0g_x=0$.
We remove the period average only from the quadratic correction,
where the constant in norm equivalence is harmless. The leading
term in the bound retains the exact Gaussian multiplier,
\begin{equation}\label{eq:outgoing-value}
 \begin{aligned}
 g'_x
 &=M^2\mathbf E_0g_x
       +(\id-\Pi)(G_x-M^2\mathbf E_0g_x),\\
 \|g'_x\|_*
 &\le M^2\|\mathbf E_0g_x\|_*+C_M\zeta_\ell^2\\
 &\le M^2e^{-\frac\beta2(C_{t+T}-C_t)_{xx}}
          [2\ell^2|z_\ell|+CB\zeta_\ell^2]
       +C_M\zeta_\ell^2\\
 &\le2\ell'^2|z_{\ell'}|+[CM^\alpha B+C_M]\zeta_\ell^2.
 \end{aligned}
\end{equation}
Indeed, the linear term agrees with the Gaussian coefficient
at the next scale,
\[
 M^2e^{-\frac\beta2(C_{t+T}-C_t)_{xx}}\cdot2\ell^2|z_\ell|
 =2\ell'^2|z_{\ell'}|.
\]
For the quadratic term we use
$e^{-\beta(C_{t+T}-C_t)_{xx}/2}\le C_\beta M^{-\beta/(4\pi)}$.

For the derivative at $x$, period averaging commutes with
derivatives and contracts the supremum norm. The density estimate gives
\begin{equation}\label{ap:central-first}
 \begin{aligned}
 \|\partial_xg'_x\|_\infty
 &\le2\|\partial_xG_x\|_\infty
 \le2\sqrt2M^2\|g_x\|_\infty
               \sqrt{(Q_Te_x,C_T^{-1}Q_Te_x)}\\
 &\le C_\beta M^2\zeta_\ell\ell'^{-2}
 =C_\beta M^{2-\alpha}\zeta_{\ell'}\ell'^{-2}.
 \end{aligned}
\end{equation}
This bound is independent of the previous $K,H,B$ because
$\|g_x\|_\infty\le(2C_\beta+1)\zeta_\ell$.

For the full Hessian, differentiate each scalar $G_x$ before summing
over $x$. This retains the distance to $x$ in both derivative indices.
With $F=N+\lambda R$, the conditional derivative formula gives
\[
 \partial_j\partial_i\mathbf E_\lambda g_x
 =\partial_j\mathbf E_\lambda[Q_T\nabla g_x]_i
       -\partial_j\operatorname{Cov}_\lambda(g_x,[Q_T\nabla F]_i).
\]
The first observable uses only one derivative of $g_x$,
\[
 \begin{aligned}
 \|[Q_T\nabla g_x]_i\|_\infty
 &\le\sum_aQ_T(i,a)\|\partial_ag_x\|_\infty\\
 &\le C(K+1)\zeta_\ell\ell^{-2}
          \sum_aQ_T(i,a)e^{-6\mu d(a,x)/\ell}\\
 &\le C(K+1)\zeta_\ell\ell'^{-2}e^{-6\mu d(i,x)/\ell'}.
 \end{aligned}
\]
We put its next derivative on the conditional density. For the covariance
we use~\eqref{fc:force-covariance-derivative}. Thus
\[
 \begin{aligned}
 |\partial_j\mathbf E_\lambda[Q_T\nabla g_x]_i|
 &\le\sqrt{2(Q_Te_j,C_T^{-1}Q_Te_j)}
                   \|[Q_T\nabla g_x]_i\|_\infty\\
 &\le C(K+1)\zeta_\ell\ell'^{-4}e^{-6\mu d(i,x)/\ell'},\\
 |\partial_j\operatorname{Cov}_\lambda(g_x,[Q_T\nabla F]_i)|
 &\le C_{M,K}H\zeta_\ell^2
                 \ell'^{-4}e^{-6\mu d(i,x)/\ell'}.
 \end{aligned}
\]
Since $H\zeta_\ell\le1$, integration in $\lambda$ gives
$|\partial_j\partial_iG_x|\le C_{M,K}\zeta_\ell\ell'^{-4}
e^{-6\mu d(i,x)/\ell'}$, with $M^2$ absorbed in $C_{M,K}$.
Gaussian convolution makes $G_x$ smooth in its background. Symmetry
of this scalar Hessian gives the same bound with $i,j$ interchanged,
and hence decay in both indices,
\[
 \begin{aligned}
 |\partial_i\partial_jG_x|
 &\le C_{M,K}\zeta_\ell\ell'^{-4}
      \min\{e^{-6\mu d(i,x)/\ell'},e^{-6\mu d(j,x)/\ell'}\}\\
 &\le C_{M,K}\zeta_\ell\ell'^{-4}
                 e^{-3\mu[d(i,x)+d(j,x)]/\ell'}.
 \end{aligned}
\]
We can now sum directly with the weight needed for the next step.
The triangle inequality gives
\[
 \begin{aligned}
 \ell'^{-2}\sum_{x,j}e^{\mu d(i,j)/(2\ell')}
                  \sup_{\varphi'}|\partial_i\partial_jG_x|
 &\le C_{M,K}\zeta_\ell\ell'^{-6}
       \sum_x e^{-5\mu d(i,x)/(2\ell')}
       \sum_j e^{-5\mu d(j,x)/(2\ell')}\\
 &\le C_{M,K}\zeta_\ell\ell'^{-2}.
 \end{aligned}
\]
Before projection, the full potential is the neutral convolution
plus $\ell'^{-2}\sum_xG_x$. Combining the last bound with
\eqref{ap:neutral-recursion} therefore proves
\begin{equation}\label{ap:full-hessian-reset}
 \begin{aligned}
 &\sup_i\sum_j e^{\mu d(i,j)/(2\ell')}
       \sup_{\varphi'}|\partial_i\partial_j(N'+\ell'^{-2}\sum_xg'_x)|\\
 &\le[C_\mu H+C_{M,K}+C_{\beta,\mu,M}H^2\zeta_\ell]
                    \zeta_\ell\ell'^{-2}.
 \end{aligned}
\end{equation}

\subsection{First derivatives of the local terms}
To estimate the phase-fixed derivative of $g'_x$, we apply
\eqref{ap:gaussian-phase-error} to $g_x$ over the first half of the
time step. Applying it to $\partial_i^{(x)}g_x$ instead would require
second derivatives of $g_x$. We differentiate the remaining Gaussian
density to estimate the derivative of the convolved error.

\begin{lemma}\label{lem:local-periodic-derivatives}
Assume $M\ge\max\{256,192/\alpha\}$,
\eqref{ap:full-hessian-row}--\eqref{ap:norms}, and the scale smallness
assumptions. Then, uniformly in field, cutoff and volume,
\begin{equation}\label{ap:periodic-derivatives}
 \|\partial_i^{(x)}(G_x-\Pi G_x)\|_*
 \le[CM^{3\alpha/2}B+C_M(1+B^{\alpha/4})]\zeta_\ell^2\ell'^{-2}
                       e^{-6\mu d(i,x)/\ell'}.
\end{equation}
The constant $C$ is independent of $M,K,H$.
\end{lemma}
\begin{proof}
We first compensate a change at $i$ by a change at $x$. This keeps
the phase fixed using a localized direction. After estimating it,
\eqref{fc:local-direction-conversion} will give the derivative
$\partial_i^{(x)}$, whose compensation uses a common field shift.
Under the convolution the localized direction is
\[
 \begin{aligned}
 w_i&=Q_Te_i-\frac{Q_{t+T}(i,x)}{Q_{t+T}(x,x)}Q_Te_x,\\
 (Q_te_x,w_i)
 &=Q_{t+T}(i,x)
       -\frac{Q_{t+T}(i,x)}{Q_{t+T}(x,x)}Q_{t+T}(x,x)=0.
 \end{aligned}
\]
Thus $\partial_{w_i}g_x=\sum_a(w_i)_a\partial_a^{(x)}g_x$,
where the phase-fixed derivatives on the right use the incoming time $t$.
With $F=N+\lambda R$, differentiation of the normalized conditional
average gives
\begin{equation}\label{fc:conditional-first}
 \left(\partial_i-
       \frac{Q_{t+T}(i,x)}{Q_{t+T}(x,x)}\partial_x\right)
       \mathbf E_\lambda g_x
 =\mathbf E_\lambda\partial_{w_i}g_x
          -\operatorname{Cov}_\lambda(g_x,\partial_{w_i}F).
\end{equation}
The force covariance was estimated in~\eqref{fc:force-covariance}.
The heat-kernel upper bound and diagonal lower bound give
$Q_{t+T}(i,x)/Q_{t+T}(x,x)\le Ce^{-6\mu d(i,x)/\ell'}$, so
\[
 \begin{aligned}
 |\operatorname{Cov}_\lambda(g_x,\partial_{w_i}F)|
 &\le|\operatorname{Cov}_\lambda(g_x,[Q_T\nabla F]_i)|
   +\frac{Q_{t+T}(i,x)}{Q_{t+T}(x,x)}
                |\operatorname{Cov}_\lambda(g_x,[Q_T\nabla F]_x)|\\
 &\le C_M\zeta_\ell^2\ell'^{-2}e^{-6\mu d(i,x)/\ell'}.
 \end{aligned}
\]
Only the first term of~\eqref{fc:conditional-first} still needs damping.

The incoming derivative bound gives stronger spatial decay than
the conclusion requires. The heat-kernel estimates apply with rate
$24\mu/\alpha\le\rho_0$, and $M\ge192/\alpha$ ensures
$24\mu/(\alpha\ell')\le\mu/(8\ell)$, so this weight is absorbed
by the incoming decay. Consequently,
\begin{equation}\label{fc:local-direction-bounds}
 \begin{aligned}
 \|\partial_{w_i}g_x\|_*
 &\le B\zeta_\ell^2\ell^{-2}\sum_a
       \left[Q_T(a,i)+\frac{Q_{t+T}(i,x)}{Q_{t+T}(x,x)}Q_T(a,x)\right]
                               e^{-6\mu d(a,x)/\ell}\\
 &\le CB\zeta_\ell^2\ell^{-2}\ell'^{-2}
       e^{-24\mu d(i,x)/(\alpha\ell')}
                               \sum_a e^{-3\mu d(a,x)/\ell}\\
 &\le CB\zeta_\ell^2\ell'^{-2}e^{-24\mu d(i,x)/(\alpha\ell')}.
 \end{aligned}
\end{equation}

Put $S=T/2$ and apply~\eqref{ap:gaussian-phase-error} to $g_x$,
with the comparison time $T$ replaced by $S$. Let $f_0$ be the
initial endpoint of that interpolation, where the phase is averaged
independently. The incoming bound gives
\[
 \|\mathbf E_{C_S}g_x-f_0\|_*
 \le CB\zeta_\ell^2e^{-\beta(C_{t+S}-C_t)_{xx}/2}.
\]
The remaining covariance $C_T-C_S$ lets us differentiate this error
through the Gaussian density. By~\eqref{eq:fresh-gaussian-cost},
\[
 \begin{aligned}
 (w_i,(C_T-C_S)^{-1}w_i)
 &\le2(Q_Te_i,(C_T-C_S)^{-1}Q_Te_i)\\
 &{}+2\frac{Q_{t+T}(i,x)^2}{Q_{t+T}(x,x)^2}
                  (Q_Te_x,(C_T-C_S)^{-1}Q_Te_x)\\
 &\le C\ell'^{-4}
      \left(1+\frac{Q_{t+T}(i,x)^2}{Q_{t+T}(x,x)^2}\right)
 \le C\ell'^{-4}.
 \end{aligned}
\]
The derivative commutes with the Gaussian and angular averages
defining $f_0$, and the independent phase retains its contraction,
\[
 \|\partial_{w_i}\mathbf E_{C_T-C_S}f_0\|_*
 \le e^{-\beta(C_{t+S}-C_t)_{xx}/2}
                                      \|\partial_{w_i}g_x\|_*.
\]
For the error, we may either put the derivative on this density or
commute it through the two averages. The two bounds give
\[
 \begin{aligned}
 \left\|\partial_{w_i}\mathbf E_{C_T-C_S}
                     (\mathbf E_{C_S}g_x-f_0)\right\|_*
 &\le\min\left\{
     \sqrt{(w_i,(C_T-C_S)^{-1}w_i)}
          \|\mathbf E_{C_S}g_x-f_0\|_*,
          2\|\partial_{w_i}g_x\|_*\right\}\\
 &\le CB\zeta_\ell^2\ell'^{-2}
       \min\{e^{-\beta(C_{t+S}-C_t)_{xx}/2},
                         e^{-24\mu d(i,x)/(\alpha\ell')}\}.
 \end{aligned}
\]
The inequality
$\min\{a,b\}\le a^{1-\alpha/4}b^{\alpha/4}$ retains the spatial rate
$6\mu$ while using only a fraction of the damping. Since
$\beta(C_{t+T}-C_{t+S})_{xx}\le C_\beta$, the covariance logarithm gives
\[
 \|\mathbf E_{C_T}\partial_{w_i}g_x\|_*
 \le CB\zeta_\ell^2\ell'^{-2}M^{-\beta(1-\alpha/4)/(4\pi)}
                         e^{-6\mu d(i,x)/\ell'}.
\]

Apply~\eqref{cr:centered-gaussian-first} to compare the conditional
average with a Gaussian average centred at the conditional mean.
The latter satisfies the bound just proved, which is uniform in the background.
The inverse-covariance estimate also holds with $C_T^{-1}$ because
$C_T-C_S\le C_T$. Hence
\[
 \begin{aligned}
 \left|\mathbf E_\lambda\partial_{w_i}g_x
       -\mathbf E_{C_T}\partial_{w_i}g_x(\mathbf E_\lambda\psi+\zeta)\right|
 &\le
 CT\|\Hess F\|_{\infty,\mathrm{op}}\ell'^{-2}
                       \sum_k\|\partial_k g_x\|_\infty\\
 &\le C_{M,K}H\zeta_\ell^2\ell'^{-2}\\
 &\le C_M\zeta_\ell^2\ell'^{-2}.
 \end{aligned}
\]
Here the constant is independent of $B$. Contraction of each average
and~\eqref{fc:local-direction-bounds} give a bound proportional to
$B$, with spatial decay. Interpolation keeps only the power
$B^{\alpha/4}$,
\[
 \begin{aligned}
 \left|\mathbf E_\lambda\partial_{w_i}g_x
       -\mathbf E_{C_T}\partial_{w_i}g_x(\mathbf E_\lambda\psi+\zeta)\right|
 &\le C_M\zeta_\ell^2\ell'^{-2}
          \min\{1,B e^{-24\mu d(i,x)/(\alpha\ell')}\}\\
 &\le C_MB^{\alpha/4}\zeta_\ell^2\ell'^{-2}e^{-6\mu d(i,x)/\ell'}.
 \end{aligned}
\]
Substitute this comparison, the Gaussian estimate and the force
covariance bound into~\eqref{fc:conditional-first} to obtain
\[
 \begin{aligned}
 &\left\|\left(\partial_i-
       \frac{Q_{t+T}(i,x)}{Q_{t+T}(x,x)}\partial_x\right)
                   \mathbf E_\lambda g_x\right\|_\infty\\
 &\le\left[CBM^{-\beta(1-\alpha/4)/(4\pi)}
              +C_M(1+B^{\alpha/4})\right]
            \zeta_\ell^2\ell'^{-2}e^{-6\mu d(i,x)/\ell'}.
 \end{aligned}
\]

We now express the outgoing derivative $\partial_i^{(x)}$, at time
$t+T$, in terms of the directions just estimated. Since
$\sum_kQ_{t+T}(k,x)=c_{t+T}/\sqrt\beta$, the two coefficients of
$\partial_x$ cancel,
\begin{equation}\label{fc:local-direction-conversion}
 \begin{aligned}
 &\left(\partial_i-
       \frac{Q_{t+T}(i,x)}{Q_{t+T}(x,x)}\partial_x\right)
 -\frac{\sqrt\beta Q_{t+T}(i,x)}{c_{t+T}}
       \sum_k\left(\partial_k-
       \frac{Q_{t+T}(k,x)}{Q_{t+T}(x,x)}\partial_x\right)\\
 &=\partial_i-\frac{\sqrt\beta Q_{t+T}(i,x)}{c_{t+T}}
                        \sum_k\partial_k
 =\partial_i^{(x)}.
 \end{aligned}
\end{equation}
The coefficient of this sum compensates for the number of sites
within distance $\ell'$ of $x$,
\begin{equation}\label{fc:phase-sum}
 \frac{\sqrt\beta Q_{t+T}(i,x)}{c_{t+T}}
       \sum_k e^{-6\mu d(k,x)/\ell'}
 \le C\ell'^{-2}e^{-6\mu d(i,x)/\ell'}\ell'^2
 \le Ce^{-6\mu d(i,x)/\ell'}.
\end{equation}
Here $\sqrt\beta/c_{t+T}=e^{\eps^2m^2(t+T)/2}\le C$ in the
heat regime. This proves the same bound for
$\partial_i^{(x)}\mathbf E_\lambda g_x$, independently of volume.
After multiplication by $M^2$ in $G_x$, the exponent of its term
proportional to $B$ is
\[
 2-\frac\beta{4\pi}(1-\alpha/4)
 =\alpha+\frac{\alpha(2-\alpha)}4
 \le\frac{3\alpha}2<2\alpha.
\]
The term proportional to $B$ therefore grows more slowly than the
induction bound $B\zeta_{\ell'}^2=BM^{2\alpha}\zeta_\ell^2$.
Finally $\partial_i^{(x)}(\id-\Pi)G_x=(\id-\Pi)\partial_i^{(x)}G_x$.
Projection and norm equivalence act only on this quadratic error,
and give~\eqref{ap:periodic-derivatives}.
\end{proof}

\subsection{Completion of the estimates}
\label{sec:completion}

\begin{proof}[Proof of Proposition~\ref{prop:outgoing-decomposition}]
Fix $\mu$ as above. Let $C$ cover the propagated quadratic coefficients
independently of $M,K,H$, and let $C_\mu$ include the linear coefficient
in~\eqref{ap:full-hessian-reset}. Choose
\begin{equation}\label{eq:M-choice}
 M\ge\max\{256,192/\alpha\},\qquad
 M^{\alpha/2}\ge4C,\qquad M^\alpha\ge4C_\mu.
\end{equation}
First choose $K\ge\max\{1,2\sqrt\beta,C_\beta M^{2-\alpha}\}$
using~\eqref{ap:central-first}, then
$H\ge\max\{1,2\beta,4C_{M,K}M^{-\alpha}\}$
using~\eqref{ap:full-hessian-reset}.
Let $C_M\ge1$ cover all remaining errors, including the fixed $K,H$.
The choice of $M$ controls the term proportional to $B$. Since
$\alpha/4<1$, the remaining terms divided by $B$ tend to zero.
We can therefore choose $B$ large enough that
\[
 CM^{-\alpha/2}
 +\frac{C_M(1+B^{\alpha/4})}{BM^{2\alpha}}
 \le\tfrac14+\tfrac14=\tfrac12.
\]
Equivalently,
\begin{equation}\label{ap:periodic-recursion}
 [CM^{3\alpha/2}B+C_M(1+B^{\alpha/4})]\zeta_\ell^2
 \le\tfrac12B\zeta_{\ell'}^2.
\end{equation}
Finally choose $0<\zeta_*\le1$ such that
\begin{equation}\label{eq:curvature-smallness}
 B\zeta_*\le1,\qquad C_MH\zeta_*\le\tfrac18.
\end{equation}
The last choice ensures conditional convexity and bounds the
quadratic term in~\eqref{ap:full-hessian-reset}. Substituting our
choices in that estimate gives
\[
 \begin{aligned}
 [C_\mu H+C_{M,K}+C_{\beta,\mu,M}H^2\zeta_\ell]
                                  \zeta_\ell\ell'^{-2}
 &\le H(M^\alpha/4+M^\alpha/4+1/8)\zeta_\ell\ell'^{-2}\\
 &\le H\zeta_{\ell'}\ell'^{-2}.
 \end{aligned}
\]
The derivative at $x$ is covered by our choice of $K$,
independently of its previous bound. Since $M^\alpha\le M^{3\alpha/2}$,
the same choice of $B$ covers the quadratic term
in~\eqref{eq:outgoing-value}. Equation~\eqref{ap:periodic-recursion}
gives the phase-fixed derivative bound.
Reality, periodicity and zero average are preserved by the construction.
Every constant has now been fixed in terms of $\beta$, before cutoff
or volume is chosen.
\end{proof}

\begin{proof}[Proof of Proposition~\ref{prop:finite-scale}]
For $z=0$ there is nothing to estimate. At $\ell=1$, the cosine has
zero phase-fixed derivatives, norm $2|z|\eps^\alpha$, derivative at
$x$ at most $2\sqrt\beta\zeta_1$, and full Hessian row at most
$2\beta\zeta_1$. Our choices of $K,H$ cover these bounds.
Choose $c_\beta$ small enough for the heat regime and smallness
conditions. Proposition~\ref{prop:outgoing-decomposition} then gives
the decomposition at lengths $M^j$.

For $\ell\ge M$, choose $j\ge0$ with $M\le\ell/M^j<M^2$.
Iterate to $M^j$, then take the final step with ratio $\ell/M^j$.
We need only the value and Hessian bounds at this last step.
In the proofs of~\eqref{eq:outgoing-value} and~\eqref{ap:full-hessian-reset},
the heat bounds are uniform in the ratio and the remaining scale
factors are bounded powers of it. In particular, the product ensuring
conditional convexity is at most
$C(\ell/M^j)^2H\zeta_{M^j}\le CM^4H\zeta_\ell$.
A single reduction of $c_\beta$ therefore makes the final step
admissible for every ratio in $[M,M^2]$. These two estimates give
\[
 \|g_x\|_*\le C_\beta\zeta_\ell,\qquad
 \sup_i\sum_j e^{\mu d(i,j)/(2\ell)}
       \sup_\varphi|\partial_i\partial_j(N+R)|
 \le C_\beta\zeta_\ell\ell^{-2}.
\]
The period average $N$ inherits this row bound. On the sine-Gordon
trajectory it is also invariant under spatial translations. Indeed,
permuting the field values by $a\mapsto a+v$ leaves $V_0$ and the
Gaussian law unchanged and commutes with $Q_t$ and period averaging.
Consequently,
\[
 \sup_\varphi|\partial_a\partial_{a+v}N(\varphi)|
 =\sup_\varphi|\partial_0\partial_vN(\varphi)|.
\]
This equality concerns field suprema; the Hessian at a fixed field
need not be spatially translation invariant.

The neutral difference identity~\eqref{ap:neutral-difference-form}
now turns the row bound into the required gradient form. Represent
each torus displacement $v$ by a shortest lattice vector,
$|v|_1=d(0,v)$, and choose a shortest path
$0=v_0,\ldots,v_{|v|_1}=v$. Then
\[
 \begin{aligned}
 \sum_a(w_{a+v}-w_a)^2
 &\le |v|_1\sum_{k=1}^{|v|_1}\sum_a
                 (w_{a+v_k}-w_{a+v_{k-1}})^2\\
 &\le |v|_1^2\|Dw\|_2^2.
 \end{aligned}
\]
The first inequality is Cauchy--Schwarz along the path. In the second,
translation of the sum over $a$ bounds each nearest-neighbour
difference by the gradient energy. Applying this estimate before
summing over displacements gives
\begin{equation}\label{ap:neutral-form-bound}
 \begin{aligned}
 |(w,\Hess N(\varphi)w)|
 &\le\tfrac12\sum_v\sup_\psi|\partial_0\partial_vN(\psi)|
                              \sum_a(w_{a+v}-w_a)^2\\
 &\le\tfrac12\sum_v d(0,v)^2
                  \sup_\psi|\partial_0\partial_vN(\psi)|\,\|Dw\|_2^2\\
 &\le C_\mu\ell^2\sum_v e^{\mu d(0,v)/(2\ell)}
                  \sup_\psi|\partial_0\partial_vN(\psi)|\,\|Dw\|_2^2\\
 &\le C_\beta\zeta_\ell\|Dw\|_2^2.
 \end{aligned}
\end{equation}
The third line uses $r^2e^{-\mu r/(2\ell)}\le C_\mu\ell^2$.
Thus the same weighted row controls both the full Hessian and the
neutral gradient form. In physical units the three bounds are
\[
 \begin{aligned}
 \|\Hess[V_t(Q_t\cdot)]\|_{\infty,\mathrm{op}}
 &\le C_\beta\zeta_\ell\ell^{-2}
  =C_\beta|z|\eps^2(\eps\ell)^{-\beta/(4\pi)},\\
 \osc R
 &\le2|\Lambda|\ell^{-2}C_\beta\zeta_\ell
  =C_\beta|z|L^2(\eps\ell)^{-\beta/(4\pi)},\\
 |(w,\Hess Nw)|
 &\le C_\beta\zeta_\ell\|Dw\|_2^2
  =C_\beta|z|(\eps\ell)^\alpha\|Dw\|_2^2.
 \end{aligned}
\]

For $1\le\ell<M$, take $N=0$ and $R=V_t(Q_t\cdot)$.
The hypothesis of~\eqref{eq:hessian-convolution} follows directly
from the bare Hessian,
\[
 \|C_t\|_{\mathrm{op}}\|\Hess V_0\|_{\infty,\mathrm{op}}
 \le 2\beta t|z|\eps^\alpha
 \le 2\beta M^2\zeta_1\le\tfrac18,
\]
after reducing $c_\beta$. For $t>0$, apply that estimate with
$\Sigma=C_t$ and $W_1=W_2=\id$, and use $\|Q_t\|_{\mathrm{op}}\le1$.
At $t=0$ the same conclusion is the bare Hessian bound.
Also $\inf V_0\le-\log\mathbf E_{C_t}e^{-V_0(\varphi+\zeta)}
\le\sup V_0$. Hence
\[
 \begin{aligned}
 \|\Hess[V_t(Q_t\cdot)]\|_{\infty,\mathrm{op}}
 &\le4\beta|z|\eps^\alpha
 \le C_\beta|z|\eps^2(\eps\ell)^{-\beta/(4\pi)},\\
 \osc R&\le\osc V_0\le4|z|L^2\eps^{-\beta/(4\pi)}
 \le C_\beta|z|L^2(\eps\ell)^{-\beta/(4\pi)}.
 \end{aligned}
\]
The factors lost here are at most powers of the fixed $M$.
The neutral form bound is zero. This proves the proposition for
every real length $\ell\ge1$.
\end{proof}

\appendix

\section{Heat-kernel estimates}\label{cr:heat-section}

The heat estimates of~\cite[Lemmas~A.1 and 3.3]{BB2021} give the following bounds with
spatial weights and differences, uniformly before the diffusion
length reaches the mass length or the side of the torus. Write $D^k$
for any $k$ nearest-neighbour differences at either endpoint, and $d$
for periodic nearest-neighbour distance. Fix $0<c_0\le1$ and assume
\begin{equation}\label{cr:heat-regime}
 \ell'\ge1,\qquad \ell'\le c_0L/\eps,\qquad \eps m\ell'\le c_0.
\end{equation}

\begin{lemma}\label{cr:heat-bounds}\label{cr:integrated-heat}\label{cr:covariance-logarithm}
There is $\rho_0>0$ such that, for $0<\rho\le\rho_0$,
$0\le s\le2\ell'^2$ and $k=0,1,2$,
\begin{align}
 \sup_x\sum_y e^{\rho d(x,y)/\ell'}|D^ke^{-sA}(x,y)|
 &\le C_\rho(1+s)^{-k/2},\label{cr:heat-schur}\\
 |D^ke^{-sA}(x,y)|
 &\le C_\rho(1+s)^{-1-k/2}e^{-\rho d(x,y)/\ell'}.
 \label{cr:reserve-point}
\end{align}
The row bound also holds for columns, and its time integral gives
\begin{equation}\label{cr:integrated-operators}
 \sup_x\sum_y e^{\rho d(x,y)/\ell'}C_r(x,y)\le C_\rho r,
 \qquad 0\le r\le2\ell'^2.
\end{equation}
For $0\le t\le\ell'^2$,
\begin{equation}\label{cr:diagonal-log}
 C_t(0,0)=\frac1{4\pi}\log(1+t)+O(1),\qquad
 (e^{-tA})_{xx}\ge\frac{c}{1+t}.
\end{equation}
All constants are independent of cutoff and volume. If
$\ell'=M\ell$, $\ell\ge1$, $M\ge2$, and $T=\ell'^2-\ell^2$, then $Q_T,DQ_T$
satisfy pointwise bounds with prefactors $C_\rho\ell'^{-2},C_\rho\ell'^{-3}$
and the same exponential decay. Also,
$C_{\ell'^2-1}(0,0)-C_{\ell^2-1}(0,0)
=(2\pi)^{-1}\log M+O(1)$, with constants independent of $M$.
\end{lemma}
\begin{proof}
The lattice heat-kernel estimate and torus periodization in
\cite[Lemma~A.1]{BB2021} give, for $p_s(x)=(e^{s\Delta_{\Z^2}})(0,x)$
and sufficiently small $\rho_0$,
\[
 e^{2\rho|x|_1/\ell'}|D^kp_s(x)|
 \le C(1+s)^{-1-k/2}e^{-c|x|_1/(2\sqrt{1+s})}.
\]
For $s<1$, the same bound follows from
$p_s\le p_1(0)^{-1}p_{s+1}$ and the expression of $D^k$ as finitely
many translates. The exponential sum is at most $C(1+s)$, which
proves the weighted row bound on $\Z^2$. On the torus,
\[
 (e^{-sA})(x,y)=e^{-\eps^2m^2s}
       \sum_{v\in(L/\eps)\Z^2}p_s(y-x+v).
\]
The periodic distance $d(x,y)$ is at most $|y-x+v|_1$ for every
lift. Summing over $y$ and its lifts therefore bounds the weighted
torus row by the infinite-lattice row. For the point estimate, retain
one rate $\rho$ for $d(x,y)$ and sum the remaining exponential over
lifts. Their separation is $L/\eps\ge\ell'/c_0$, so this sum is
uniformly bounded. The mass factor is at most one. Time integration gives~\eqref{cr:integrated-operators},
and $T/2\asymp\ell'^2$ gives the stated bounds on $Q_T$ and $DQ_T$.

The covariance logarithm follows from~\cite[Lemma~3.3]{BB2021}.
Its condition~(3.17) allows either $Lm\ge1$ or
$t\le\eps^{-2}(m^{-2}\wedge L^2)$. The latter follows
from~\eqref{cr:heat-regime}. In this range its length is
$1\vee\sqrt t$, whose logarithm is $\frac12\log(1+t)+O(1)$.
Subtracting at the two times gives the increment with an $O(1)$
error independent of $M$. The lower diagonal follows from the diagonal asymptotic
in~\cite[Lemma~A.1]{BB2021}, the bound $p_t(0)\ge e^{-4t}$
on bounded intervals, positivity of periodization and
$e^{-\eps^2m^2t}\ge e^{-c_0^2}$.
\end{proof}

\section*{Acknowledgements}
The author was partially supported by a grant of the Ministry of
Research, Innovation and Digitization, CCCDI - UEFISCDI, project
number ROSUA-2024-0001, within PNCDI IV.

\section*{Declaration of generative AI use}

This paper was written by the author and was not generated by artificial
intelligence. OpenAI Codex 5.6 Sol  was used to assist the author in a first pass check for correctness of mathematical calculations from handwritten notes, as well as literature checks. All mathematical ideas, arguments, and proofs are those of
the author, who takes full responsibility for the contents of the paper.

\section*{Data availability statement}

No data were used for the research described in the article.

\section*{Conflict of interest statement}

The author has no conflict of interest to declare.

\end{document}